\documentclass{elsart}
\usepackage{graphicx,psfrag}
\usepackage{amsmath}
\usepackage{amssymb}
\begin{document}
\begin{frontmatter}

\title{Finite strain viscoplasticity with
nonlinear kinematic hardening:
phenomenological modeling and
time integration}
\author{A.V. Shutov\corauthref{cor}},
\corauth[cor]{Corresponding author. Tel.: +49-0-371-531-35024; fax: +49-0-371-531-23419.}
\ead{alexey.shutov@mb.tu-chemnitz.de}
\author{R. Krei{\ss}ig}

\address{Institute of Mechanics,
Chemnitz University of Technology,
Str. d. Nationen 62, D-09111 Chemnitz, Germany}

\begin{abstract}
This article deals with a viscoplastic material
model of overstress type.
The model is based on a
multiplicative decomposition of the deformation gradient
into elastic and inelastic part. An additional multiplicative decomposition
of inelastic part is used to describe a nonlinear kinematic hardening of Armstrong-Frederick type.

Two implicit time-stepping methods are adopted for numerical integration of
evolution equations, such that the plastic incompressibility constraint is exactly satisfied.
The first method is based on the tensor exponential.
The second method is a modified Euler-Backward method.
Special numerical tests show that both
approaches yield similar results even for finite inelastic increments.

The basic features of the material response, predicted by the material model,
are illustrated with a series of numerical simulations.
\end{abstract}
\end{frontmatter}


\section*{Notation}

\begin{tabbing}
$\mathbf F$ \quad \quad \quad \quad \quad  \quad \quad  \= deformation gradient \\
$\mathbf F_{\text{i}}$ \>  inelastic part of the deformation gradient \\
$\hat{\mathbf F}_{\text{e}}$ \>   elastic part of the deformation gradient \\
$\mathbf F_{\text{ii}}$ \> dissipative part of  $\mathbf F_{\text{i}}$ \\
$\check{\mathbf F}_{\text{ie}}$ \> energy storage part of $\mathbf F_{\text{i}}$ \\
$\mathcal{K}$ \> current configuration \\
$\tilde{\mathcal{K}}$ \>  reference configuration \\
$\hat{\mathcal{K}}$ \> stress-free intermediate configuration \\
$\check{\mathcal{K}}$\> intermediate configuration of microstructure \\
$\mathbf C$ \> right Cauchy-Green tensor (see \eqref{reicg}) \\
$\mathbf C_{\text{i}}$ \> inelastic right Cauchy-Green tensor (see $\eqref{inreicg}_1$) \\
$\mathbf C_{\text{ii}}$ \> inelastic right Cauchy-Green tensor of microstructure (see $\eqref{inreicg}_2$) \\
$\hat{\mathbf C}_{\text{e}}$ \> elastic right Cauchy-Green tensor (see $\eqref{cecie}_1$) \\
$\check{\mathbf C}_{\text{ie}}$ \> elastic right Cauchy-Green tensor of microstructure (see $\eqref{cecie}_2$) \\
$\mathbf E$ \> Green strain tensor (see $\eqref{totdef}_1$) \\
$\mathbf \Gamma$ \> Almansi strain tensor (see $\eqref{totdef}_2$) \\
$\mathbf T$ \> Cauchy stress tensor \\
$\mathbf S$ \> weighted Cauchy tensor (Kirchhoff tensor) (see \eqref{Kirch}) \\
$\hat{\mathbf S}$, $\tilde{\mathbf T}$ \> 2nd Piola-Kirchhoff tensors operating on
$\hat{\mathcal{K}}$, $\tilde{\mathcal{K}}$, respectively (see \eqref{Kirch2}) \\
$\check{\mathbf X}$, $\hat{\mathbf X}$, $\tilde{\mathbf X}$\> backstress tensors
operating on $\check{\mathcal{K}}$, $\hat{\mathcal{K}}$ and $\tilde{\mathcal{K}}$, respectively (see \eqref{backst1}) \\
$\hat{\mathbf \Sigma}$ \> the driving force for inelastic flow  (see $\eqref{drforc}_1$) \\
$\check{\mathbf \Xi}$ \>  the driving force for inelastic flow of microstructure (see $\eqref{drforc}_2$)  \\
$\mathbf 1$ \> second-rank identity tensor \\
$\mathbf M^*, \mathbf M_*$ \> covariant pull-back and push-forward (see \eqref{cova})  \\
$({\mathbf M^{-\text{T}}})^*,
({\mathbf M^{-\text{T}}})_*$ \> contravariant pull-back and push-forward (see \eqref{contra})  \\
$\mathbf A \cdot \mathbf B = \mathbf A \mathbf B$ \> product (composition) of two second-rank tensors \\
$\mathbf A : \mathbf B$ \> scalar product of two second-rank tensors (see \eqref{scalprod}) \\
$\| \mathbf A \|$ \> $l_2$ norm of a second-rank tensor (Frobenius norm) (see $\eqref{defi}_1$) \\
$\| \mathbf A \|^*$ \> induced norm of a second-rank tensor (spectral norm)  (see \eqref{Openo}) \\
$(\cdot)^{\text{D}}$ \> deviatoric part of a tensor (see $\eqref{defi}_2$) \\
$ (\cdot)^{\text{T}}$  \> transposition of a tensor \\
$ (\cdot)^{-\text{T}}$ \>  inverse of transposed \\
$\text{tr}(\cdot)$ \> trace of a second-order tensor \\
$\stackrel{\triangle} {(\cdot)}$ \> covariant Oldroyd rate with respect to $\hat{\mathcal{K}}$  (see $\eqref{oldro}_1$)\\
$\stackrel{\diamondsuit} {(\cdot)}$ \> covariant Oldroyd rate with respect to $\check{\mathcal{K}}$ (see $\eqref{oldro}_2$) \\
$\overline{(\cdot)}$ \> unimodular part of a tensor (see \eqref{unim}) \\
$\text{sym}(\cdot)$ \> symmetric part of a tensor (see \eqref{sympart})  \\
$\text{skew}(\cdot)$ \> skew-symmetric part of a tensor (see $\eqref{Sympr4}_2$)  \\
$\langle x \rangle$ \> MacCauley bracket (see $\eqref{perz}_3$) \\
$\psi$ \>  specific free energy \\
$\delta_{\text{i}}$ \> specific internal dissipation (see \eqref{cld}) \\
$K$ \> initial yield stress \\
$R$ \> isotropic hardening \\
${}^{\text{t}} R$ \> trial isotropic hardening (see $\eqref{excl1}_2$)\\
$s$ \> inelastic arc length \\
$s_{\text{d}}$ \> dissipative part of $s$ \\
$s_{\text{e}}$ \> energy storage part of $s$ (see \eqref{iso1}) \\
$\lambda_{\text{i}}$ \> proportionality factor (inelastic multiplier)  (see $\eqref{perz}_1$) \\
$f$ \> overstress  (see $\eqref{perz}_2$) \\
$Sym$ \> space of symmetric second-rank tensors \\
$\mathfrak{F}$ \> norm of the driving force (see \eqref{nordiv}) \\
$\xi$ \>  incremental inelastic parameter (see \eqref{ineinc})\\
$\rho_{\scriptscriptstyle \text{R}}$ \> mass density in the reference configuration \\
$k$ \> bulk modulus (see \eqref{spec1})  \\
$\mu$ \>  shear modulus (see \eqref{spec1}) \\
$c$ \> bulk modulus of microstructure (see $\eqref{spec2}_1$) \\
$\gamma$ \> hardening modulus (see $\eqref{spec2}_2$) \\
$\varepsilon$ \> technical strain (see $\eqref{uniax}_1$)\\
$\phi$ \> shear strain (see $\eqref{tors}_1$)\\
$\sigma$, $\tau$ \> axial and shear stresses, respectively\\
\end{tabbing}

\section{Introduction}

New materials, such as ultrafine-grained-aluminium (see the papers \cite{Ha1}, \cite{Ha2}),
are of special interest for many practical applications.
To promote the innovation of the new materials,
the robust numerical simulation of the material response
is required.
It is desirable to have a phenomenological description of the material which on the one hand
takes important phenomena into account, and on the other hand enables stable
numerical computations.

In this paper we investigate the simulation of rate-dependent material
behavior with equilibrium hysteresis effect (for the general introduction to
the theory of viscoplasticity see, for example, \cite{Perzyna1}, \cite{LemaitreChaboche}, \cite{Haupt}).

The Bauschinger effect is observed in most metals under non-monotonic loading.
The most popular approach to describe the Bauschinger effect
was proposed by Armstrong and Frederick \cite{Arms} in 1966.
Application of the Armstrong-Frederick hardening concept within the framework of Perzyna type viscoplasticity
(see \cite{Perzyna1}) yields
the classical material model of overstress type (see \cite{ChabocheRousselier}, \cite{LemaitreChaboche}).
This model has the advantage that it admits simple rheological interpretation (see fig. \ref{fig1}.a).
Such phenomena as creep, relaxation and nonlinear kinematic hardening are
taken into account by the model.
Simple modification of this model is
possible to include isotropic hardening as well \footnote{The diagram
in fig. \ref{fig1}.a provides insight into the
rheological modeling of kinematic hardening.
To the best of our knowledge, there is no
simple rheological diagram of viscoplastic
material with isotropic hardening.}.

Several strategies can be adopted
for the generalization of this model to finite strains (see, for example,
\cite{DogSid}, \cite{Tsakmak}, \cite{Luhrs}, \cite{Svendsen}, \cite{Lion}, \cite{Helm1}, \cite{Mollicaetal}). Some of the
generalizations were analyzed numerically in \cite{DettRes}. Following the elegant approach of Lion \cite{Lion}, we use the rheological
interpretation (fig. \ref{fig1}.a) of the classical model to construct its finite-strain counterpart.

The specific assumptions of the material modeling used in this paper are as follows:
\begin{itemize}
\item Multiplicative decomposition of the deformation gradient
into elastic and inelastic part: $\mathbf F = \hat{\mathbf F}_{\text{e}} \mathbf F_{\text{i}}$ (\cite{Kroen}, \cite{Lee}).
\item Multiplicative decomposition of the inelastic part
into energy storage part and dissipative part: $\mathbf F_{\text{i}} = \check{\mathbf F}_{\text{ie}} \mathbf F_{\text{ii}}$ (\cite{Lion}).
\item Free energy is a sum of appropriate isotropic strain energy functions (\cite{Lion}).
\end{itemize}
The resulting material model takes both kinematic and isotropic hardening into account.
The thermodynamic consistency is proved.

The purpose of the present paper is threefold. First,
we formulate the material model under consideration. In particular, we transform
the constitutive equations to the reference configuration in order to simplify
the numerical treatment.
Next, two implicit
schemes for the numerical integration of evolution equations are developed. Finally,
we analyse numerically the basic properties of the material response, predicted by the model.

A global \emph{implicit} time stepping procedure in the context of displacement based FEM
requires a
proper stress algorithm (local integration algorithm) \cite{Zienkiewic}.
Such algorithm provides the stresses and the
consistent tangent operator as a function of the strain history locally at each integration point.
A set of internal variables is used in this paper to describe the history dependence, and
the stress algorithm includes \emph{implicit} integration
of a system of differential (evolution) and algebraic equations.

Two most popular implicit schemes for integration of inelastic strains
in the context of viscoplasticity/plasticity are:
\begin{itemize}
\item \emph{Backward-Euler scheme}, also referred as \emph{implicit Euler scheme} (see, for example,
\cite {Hartmann}, \cite{SimHug}, \cite{SimMi}, \cite{Helm1}, \cite{DettRes}).
\item \emph{Exponential scheme}, also referred as \emph{Euler scheme with exponential map}
(see, for example, \cite{WebAnan}, \cite{MiStei}, \cite{Miehe}, \cite{DettRes}).
\end{itemize}
The exponential scheme is advantageous since it retains the inelastic
incompressibility even for finite time steps. Thus, an important geometric property
of the solution is automatically preserved.
Moreover, the numerical error of Euler-Backward method, related to the
violation of incompressibility, \emph{tends to accumulate over time} (see, for example, \cite{DettRes}, \cite{Helm2}).
Therefore, even for small time steps,
the numerical solution deviates from the exact solution after some period of time.

Helm \cite{Helm2} modified the classical
Euler-Backward scheme, using a projection on the group of unimodular tensors, to enforce
the incompressibility of inelastic flow.

In this work we implement in a uniform manner both modified Euler-Backward method (MEBM) and the
exponential method (EM). Both methods result in a nonlinear system of equations
with respect to strain-like internal variables $\mathbf{C}_{\text{i}} = \mathbf{F}_{\text{i}}^{\text{T}} \mathbf{F}_{\text{i}}$,
$\mathbf{C}_{\text{ii}} = \mathbf{F}_{\text{ii}}^{\text{T}} \mathbf{F}_{\text{ii}}$ and  $\xi= \lambda_{\text{i}} \ \Delta t$
\footnote{$\xi \geq 0$ is an incremental inelastic parameter, defined by \eqref{ineinc}}.
This nonlinear system is
split into two subproblems: \\
First subproblem: Finding $\mathbf{C}_{\text{i}}, \mathbf{C}_{\text{ii}}$ with a given $\xi$. \\
Second subproblem: Finding $\xi$, such that an incremental consistency condition is satisfied. \\
This adapted strategy is
more robust than the straightforward application of a nonlinear solver to the original
system of equations. At the same time, this approach is
not limited by the special form of the free energy, and finite elastic strains
are likewise allowed.
Moreover, the stress algorithms are applicable in the limiting case of
rate-independent plasticity (as viscosity tends to zero).

Although the material response is anisotropic, it is shown that MEBM as well as EM
\emph{exactly preserve the symmetry} of $\mathbf{C}_{\text{i}}$ and $\mathbf{C}_{\text{ii}}$.
Furthermore, the accuracy of both integration algorithms is verified
with the help of special numerical tests. Both methods provide similar results with almost
the same integration error. A common feature of MEBM and EM is that the numerical error
is not accumulated over time.

The phenomenological description of each specific material
can be schematically subdivided into three steps:
\begin{itemize}
\item Material testing, such that the important phenomena make themselves evident.
\item Choosing an appropriate phenomenological model, that reproduces qualitatively the experimental data.
\item Parameter identification, using the experimental data.
\end{itemize}

To illustrate
the basic characteristics of the material model we simulate a
series of material testing experiments. These experiments are uniaxial tension and
torsion under monotonic and cyclic loading. In particular, we conclude that
the material model can be used (after a proper parameter identification) to describe the
mechanical response of an aluminium alloy processed by ECA-pressing \cite{Ha1}, \cite{Ha2}.

Throughout this article, bold-faced symbols denote first- and second-rank tensors in $\mathbb{R}^3$.
Expression $a := b$ means $a$ is defined to be another name for $b$.

\section{Material model of finite viscoplasticity}

The material model is motivated by the rheological diagram in fig. \ref{fig1}.a.
This diagram takes
the kinematic hardening of Armstrong-Frederick type into account
(for the sake of simplicity the isotropic hardening is omitted in the diagram).
The
total inelastic strains and the inelastic strains of
microstructure are used as internal variables.
The evolution of these quantities is closely related to
the energy dissipation during the inelastic processes.
Besides, additional real-valued strain-like internal variables are introduced in order
to describe a nonlinear isotropic hardening.

\subsection{Kinematics}

\begin{figure}\centering
\psfrag{A}[m][][1][0]{\footnotesize energy}
\psfrag{B}[m][][1][0]{\footnotesize storage}
\psfrag{C}[m][][1][0]{\footnotesize dissipation}
\psfrag{D}[m][][1][0]{microstructure}
\psfrag{I}[m][][1][0]{$\tilde{\mathcal{K}}$}
\psfrag{J}[m][][1][0]{$\check{\mathcal{K}}$}
\psfrag{K}[m][][1][0]{$\hat{\mathcal{K}}$}
\psfrag{S}[m][][1][0]{$\mathcal{K}$}
\psfrag{L}[m][][1][0]{$\mathbf F_{\text{ii}}$}
\psfrag{M}[m][][1][0]{$\mathbf F_{\text{i}}$}
\psfrag{N}[m][][1][0]{$\check{\mathbf F}_{\text{ie}}$}
\psfrag{O}[m][][1][0]{$\hat{\mathbf F}_{\text{e}}$}
\psfrag{P}[m][][1][0]{$\mathbf F$}
\psfrag{Q}[m][][1][0]{a}
\psfrag{R}[m][][1][0]{b}
\scalebox{0.95}{\includegraphics{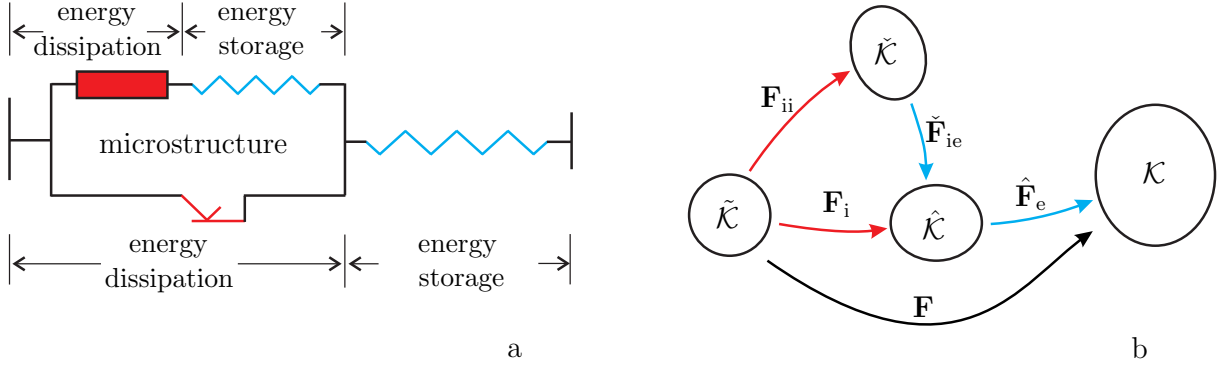}}
\caption{Rheological model (a) and a commutative diagram (b) showing corresponding
configurations with transformations of material line elements. \label{fig1}}
\end{figure}

For a fixed time instant $t \geq 0$ let $\mathcal{K} \subset \mathbb{R}^3$ be a current configuration occupied by
the solid. Suppose $\tilde{\mathcal{K}} \subset \mathbb{R}^3$ is the reference configuration,
which uniquely designates the material points. Let us consider the motion law in the form
$\mathbf \chi_t : \tilde{\mathcal{K}} \rightarrow \mathcal{K}$. For every point
$\mathbf P \in \tilde{\mathcal{K}}$ we define
the deformation gradient tensor
$\mathbf F := \frac{\displaystyle \partial \mathbf \chi_t (\mathbf P)}{\displaystyle \partial \mathbf P}$.
The deformation gradient $\mathbf F$ transforms a material line element $d \mathbf X$ on the
reference configuration $\widetilde{\mathcal{K}}$ into a current material line element $d \mathbf x$
\begin{equation}\label{tran}
d \mathbf x = \mathbf F  \ d \mathbf X.
\end{equation}
Let us consider the classical multiplicative decomposition of the deformation gradient $\mathbf F$
into \emph{elastic part} $\hat{\mathbf F}_{\text{e}}$ and \emph{inelastic part} $\mathbf F_{\text{i}}$ (\cite{Kroen}, \cite{Lee})
\begin{equation}\label{mude1}
\mathbf F = \hat{\mathbf F}_{\text{e}} \mathbf F_{\text{i}}.
\end{equation}
The mechanical justification uses the idea of the local
(within a neighborhood of the material point) elastic unloading.
The transformation rule \eqref{tran} is represented as a combination of two linear operators
\begin{equation*}\label{combinat}
d \mathbf x = \hat{\mathbf F}_{\text{e}}  \big(\mathbf F_{\text{i}}  \ d \mathbf X \big).
\end{equation*}
Therefore, we can interpret $\mathbf F_{\text{i}}  \ d \mathbf X$ as a fictitious material line element on
some intermediate configuration $\hat{\mathcal{K}}$ (see fig. \ref{fig1}.b). We will call this configuration
the \emph{stress-free intermediate configuration}.

A second multiplicative decomposition is introduced in order to simulate
a nonlinear kinematic hardening of Armstrong-Frederick type.
Following Lion \cite{Lion}, we decompose the
inelastic part $\mathbf F_{\text{i}}$ into \emph{energy storage
part} $\check{\mathbf F}_{\text{ie}}$ and \emph{dissipative part} $\mathbf F_{\text{ii}}$
\begin{equation}\label{mude2}
\mathbf F_{\text{i}} = \check{\mathbf F}_{\text{ie}} \mathbf F_{\text{ii}}.
\end{equation}
The energy storage part $\check{\mathbf F}_{\text{ie}}$ describes the heterogeneity
of elastic strains associated with the energy
storage on the microscale. The dissipative part $\mathbf F_{\text{ii}}$
can be attributed to slip processes on the microscale (see \cite{Lion}, \cite{Helm1} for details).
Decomposition \eqref{mude2} implements
the \emph{intermediate configuration of microstructure}\footnote
{In \cite{Lion} the similar configuration is called the \emph{intermediate configuration of
kinematic hardening}.} $\check{\mathcal{K}}$ (see fig. \ref{fig1}.b).
The commutative diagram in fig. \ref{fig1}.b summarizes both multiplicative decompositions.

In this paper we deal with second-order tensors, which
operate on configurations $\mathcal{K}$, $\tilde{\mathcal{K}}$,
$\hat{\mathcal{K}}$, $\check{\mathcal{K}}$.
Pull-back (push-forward) operations describe the transformation of tensor fields during the
change of configurations. Let
$\mathbf M \in \big\{\mathbf F, \mathbf F_{\text{i}}, \mathbf F_{\text{ii}}, \hat{\mathbf F}_{\text{e}}, \check{\mathbf F}_{\text{ie}} \big\}$ be
a linear transformation of material line elements on two different configurations.
We define corresponding \emph{pull-back and push-forward of covariant tensor field} by
\begin{equation}\label{cova}
\mathbf M^* (\cdot) : = \mathbf M^{\text{T}} (\cdot) \mathbf M, \quad
\mathbf M_* (\cdot) : = \mathbf M^{- \text{T}} (\cdot) \mathbf M^{-1}.
\end{equation}
\emph{Pull-back and push-forward of contravariant tensor field} are given by
\begin{equation}\label{contra}
({\mathbf M^{-\text{T}}})^* (\cdot)  = \mathbf M^{-1} (\cdot) \mathbf M^{-\text{T}}, \quad
({\mathbf M^{-\text{T}}})_* (\cdot)  = \mathbf M (\cdot) \mathbf M^{\text{T}}.
\end{equation}
Thus, the right Cauchy-Green tensor is a covariant pull-back of
$\mathbf 1$
\begin{equation}\label{reicg}
\mathbf C  := \mathbf F^* \mathbf 1 = \mathbf F^{\text{T}} \mathbf F.
\end{equation}
In the same manner, we define the inelastic right Cauchy-Green tensor $\mathbf C_{\text{i}}$ and
inelastic right Cauchy-Green tensor $\mathbf C_{\text{ii}}$ of microstructure:
\begin{equation}\label{inreicg}
\mathbf C_{\text{i}}  := \mathbf{F}_{\text{i}}^* \mathbf 1 = \mathbf F_{\text{i}}^{\text{T}} \mathbf F_{\text{i}}, \quad
\mathbf C_{\text{ii}}  := \mathbf {F}_{\text{ii}}^* \mathbf 1 = \mathbf F_{\text{ii}}^{\text{T}} \mathbf F_{\text{ii}}.
\end{equation}
Further, the elastic right Cauchy-Green tensor $\hat{\mathbf C}_{\text{e}}$ and the elastic
right Cauchy-Green tensor $\check{\mathbf C}_{\text{ie}}$ of microstructure
are defined by
\begin{equation}\label{cecie}
\hat{\mathbf C}_{\text{e}}  := {\hat{\mathbf F}_{\text{e}}}^* \mathbf 1 =
\hat{\mathbf F}_{\text{e}}^{\text{T}} \hat{\mathbf F}_{\text{e}}, \quad
\check{\mathbf C}_{\text{ie}}  := {\check{\mathbf F}_{\text{ie}}}^* \mathbf 1 =
\check{\mathbf F}_{\text{ie}}^{\text{T}} \check{\mathbf F}_{\text{ie}}.
\end{equation}

The tensors
\begin{equation}\label{totdef}
\mathbf E : = \frac{\displaystyle 1}{\displaystyle 2} (\mathbf C - \mathbf 1), \quad
\mathbf \Gamma := \mathbf {F}_*\mathbf E = \mathbf F^{-\text{T}} \mathbf E \mathbf F^{-1}
\end{equation}
are termed the Green strain tensor and the Almansi strain tensor, respectively.
Basing on $\mathbf E$, we define a corresponding strain tensor $\hat{\mathbf \Gamma}$, which operates on $\hat{\mathcal{K}}$
\begin{equation*}
\hat{\mathbf \Gamma} := (\mathbf F_{\text{i}})_*\mathbf E = \mathbf {F}_{\text{i}}^{-\text{T}} \mathbf E \mathbf {F}_{\text{i}}^{-1}.
\end{equation*}
Multiplicative decomposition \eqref{mude1} implements the additive decomposition of
$\hat{\mathbf \Gamma}$:
\begin{equation}\label{addco1}
\mathbf F = \hat{\mathbf F}_{\text{e}} \mathbf F_{\text{i}} \quad \Rightarrow \quad
\hat{\mathbf \Gamma} = \hat{\mathbf \Gamma}_{\text{i}} + \hat{\mathbf \Gamma}_{\text{e}},
\end{equation}
where $\hat{\mathbf \Gamma}_{\text{i}}$ is a purely inelastic Almansi tensor
\begin{equation*}
\hat{\mathbf \Gamma}_{\text{i}} := \frac{\displaystyle 1}{\displaystyle 2}
(\mathbf 1 - \mathbf F_{\text{i}}^{-\text{T}} \mathbf F_{\text{i}}^{-1})
\end{equation*}
and
$\hat{\mathbf \Gamma}_{\text{e}}$ is the elastic Green tensor
\begin{equation}\label{defgae}
\hat{\mathbf \Gamma}_{\text{e}} := \frac{\displaystyle 1}{\displaystyle 2}
(\hat{\mathbf F}_{\text{e}}^{\text{T}} \hat{\mathbf F}_{\text{e}} - \mathbf 1) =
\frac{\displaystyle 1}{\displaystyle 2}
(\hat{\mathbf C}_{\text{e}} - \mathbf 1).
\end{equation}

Analogously, multiplicative decomposition \eqref{mude2} implements the
additive decomposition of
the pull-back of the inelastic Almansi tensor $\hat{\mathbf \Gamma}_{\text{i}}$
to $\check{\mathcal{K}}$:
\begin{equation*}
\check{\mathbf \Gamma}_{\text{i}} := \check{\mathbf F}_{\text{ie}}^* \hat{\mathbf \Gamma}_{\text{i}} =
\check{\mathbf F}_{\text{ie}}^{\text{T}} \hat{\mathbf \Gamma}_{\text{i}} \check{\mathbf F}_{\text{ie}},
\end{equation*}
\begin{equation}\label{addmic}
\mathbf F_{\text{i}} = \check{\mathbf F}_{\text{ie}} \mathbf F_{\text{ii}} \quad \Rightarrow \quad
\check{\mathbf \Gamma}_{\text{i}} = \check{\mathbf \Gamma}_{\text{ii}} + \check{\mathbf \Gamma}_{\text{ie}},
\end{equation}
where
\begin{equation}\label{defgaie}
\check{\mathbf \Gamma}_{\text{ii}} := \frac{\displaystyle 1}{\displaystyle 2}
(\mathbf 1 - \mathbf F_{\text{ii}}^{-\text{T}} \mathbf F_{\text{ii}}^{-1}), \quad
\check{\mathbf \Gamma}_{\text{ie}} := \frac{\displaystyle 1}{\displaystyle 2}
(\check{\mathbf F}_{\text{ie}}^{\text{T}} \check{\mathbf F}_{\text{ie}} - \mathbf 1) =
\frac{\displaystyle 1}{\displaystyle 2}
(\check{\mathbf C}_{\text{ie}} - \mathbf 1).
\end{equation}
Finally, we  define the inelastic pull-back of $\hat{\mathbf \Gamma}_{\text{i}}$ to $\tilde{\mathcal{K}}$
\begin{equation*}
\tilde{\mathbf \Gamma}_{i} := \mathbf F_{\text{i}}^* \hat{\mathbf \Gamma}_{\text{i}}.
\end{equation*}

In this paper, the evolution of isotropic hardening is taken into account,
similar to the Armstrong-Frederick rule.
To this end, we introduce two real-valued
internal variables of strain type: $s$ and $s_{\text{d}}$.  The first variable
is the classical inelastic
arc length, and $s_{\text{d}}$ is interpreted as a dissipative part of $s$, such that
\begin{equation}\label{iso1}
s_{\text{e}}:=s - s_{\text{d}}
\end{equation}
controls the energy stored due to the isotropic hardening (see section 2.3).

\subsection{The concept of dual variables}
The formalism of dual variables developed by Haupt and Tsakmakis \cite{HaupTsa}
specifies the choice of stress and strain variables as well as their time derivatives.
According to this concept, we introduce the covariant Oldroyd rates
$\stackrel{\triangle} {(\cdot)}$, $\stackrel{\diamondsuit} {(\cdot)}$ with
respect to the stress-free configuration $\hat{\mathcal{K}}$ and the microstructural
configuration $\check{\mathcal{K}}$, respectively,
\begin{equation}\label{oldro}
\stackrel{\triangle} {(\cdot)} := {\mathbf {F}_{\text{i}}}_* \Big( \big( \mathbf {F}_{\text{i}}^* (\cdot) \big)^{\displaystyle \cdot} \Big),
\quad \quad
\stackrel{\diamondsuit} {(\cdot)} := {\mathbf {F}_{\text{ii}}}_* \Big( \big( \mathbf {F}_{\text{ii}}^* (\cdot) \big)^{\displaystyle \cdot} \Big),
\end{equation}
where $()^{\displaystyle \cdot}$ stands for material time derivative. The alternative representation
of covariant Oldroyd derivatives is as follows
\begin{equation}\label{altol}
\stackrel{\triangle} {(\cdot)} = (\cdot)^{\displaystyle {\cdot}} +
\hat{\mathbf L}^T_{\text{i}}
(\cdot) + (\cdot) \hat{\mathbf L}_{\text{i}}, \quad
\stackrel{\diamondsuit} {(\cdot)} = (\cdot)^{\displaystyle {\cdot}} +
\check{\mathbf L}^T_{\text{ii}}
(\cdot) + (\cdot) \check{\mathbf L}_{\text{ii}},
\end{equation}
\begin{equation*}\label{inelrate}
\hat{\mathbf L}_{\text{i}} := \dot{\mathbf F}_{\text{i}}
\mathbf F^{-1}_{\text{i}}, \quad \check{\mathbf L}_{\text{ii}} :=
\dot{\mathbf F}_{\text{ii}} \mathbf F^{-1}_{\text{ii}}.
\end{equation*}
Equation \eqref{oldro} yields the inelastic deformation rates
$\stackrel{\triangle}{\hat{\mathbf \Gamma}}_{\text{i}}$
and $\stackrel{\diamondsuit}{\check{\mathbf \Gamma}}_{\text{ii}}$ as symmetric parts of $\hat{\mathbf L}_{\text{i}}$
and $\check{\mathbf L}_{\text{ii}}$, respectively:
\begin{equation}\label{oldrspec}
\stackrel{\triangle}{\hat{\mathbf \Gamma}}_{\text{i}} = \text{sym}(\hat{\mathbf L}_{\text{i}}),
\quad
\stackrel{\diamondsuit}{\check{\mathbf \Gamma}}_{\text{ii}} = \text{sym}(\check{\mathbf L}_{\text{ii}}),
\end{equation}
where the symmetric part of a tensor is given by
\begin{equation}\label{sympart}
\text{sym}(\cdot):=\frac{\displaystyle 1}{\displaystyle 2} \Big( (\cdot) + (\cdot)^{\text{T}} \Big).
\end{equation}

Let $\mathbf T$ be the Cauchy stress tensor.
The weighted Cauchy tensor (or Kirchhoff stress tensor) is defined by
\begin{equation}\label{Kirch}
\mathbf S  : = (\text{det} \mathbf F) \mathbf T.
\end{equation}
Now, we define the 2nd Piola-Kirchhoff tensors operating on $\hat{\mathcal{K}}$ and
$\tilde{\mathcal{K}}$
using a contravariant pull-back of $\mathbf S$
\begin{equation}\label{Kirch2}
\hat{\mathbf S} : = \big({\hat{\mathbf F}_{\text{e}}^{-\text{T}}}\big)^* \mathbf S, \quad
\tilde{\mathbf T} : = \big({\mathbf F^{-\text{T}}}\big)^* \mathbf S =
\big({{\mathbf F}_{\text{i}}^{-\text{T}}}\big)^* \hat{\mathbf S}.
\end{equation}

The introduced stress and strain tensors form
the following conjugate pairs:
$(\mathbf S, \mathbf \Gamma)$,
$(\hat{\mathbf S}, \hat{\mathbf \Gamma})$, $(\tilde{\mathbf T}, \mathbf E)$,
such that the work and the stress power
are invariant under the change of configuration:
\begin{equation}\label{strpow1}
\mathbf S : \mathbf \Gamma = \hat{\mathbf S} : \hat{\mathbf \Gamma} =
\tilde{\mathbf T}: \mathbf E, \quad
\hat{\mathbf S} : \stackrel{\triangle} {\hat{\mathbf \Gamma}} =
\tilde{\mathbf T}: \dot{\mathbf E}.
\end{equation}
Here  "$:$ " denotes the scalar product of two second-rank tensors
\begin{equation}\label{scalprod}
\mathbf A : \mathbf B \ := \text{tr} (\mathbf A \cdot \mathbf B^{\text{T}}).
\end{equation}

Next, we denote by $\check{\mathbf X}$ the backstress tensor, which operates on the intermediate
configuration $\check{\mathcal{K}}$ of microstructure. This tensor can be interpreted as
a generalized force, associated with strain measure $\check{\mathbf \Gamma}_{\text{i}}$ and strain
rate $\stackrel{\diamondsuit} {\check{\mathbf \Gamma}}_{\text{i}}$. According to the concept of dual variables,
we define transformation rules for the backstress tensor, such that quantities
$\check{\mathbf X}: \check{\mathbf \Gamma}_{\text{i}}$ and
$\check{\mathbf X}: \ \stackrel{\diamondsuit} {\check{\mathbf \Gamma}}_{\text{i}}$
remain invariant under the change of configuration:
\begin{equation}\label{backst1}
\tilde{\mathbf X} := \big({{\mathbf F}_{\text{ii}}^{-\text{T}}}\big)^* \check{\mathbf X}=
{\mathbf F}_{\text{ii}}^{-1} \check{\mathbf X} {{\mathbf F}_{\text{ii}}^{-\text{T}}}
, \quad
\hat{\mathbf X} := \big({{\check{\mathbf F}}_{\text{ie}}^{-\text{T}}}\big)_* \check{\mathbf X}=
{\check{\mathbf F}}_{\text{ie}} \check{\mathbf X} {\check{\mathbf F}}_{\text{ie}}^{\text{T}},
\end{equation}
\begin{equation}\label{backst2}
\check{\mathbf X}: \check{\mathbf \Gamma}_{\text{i}} =
\hat{\mathbf X}: \hat{\mathbf \Gamma}_{\text{i}} =
\tilde{\mathbf X}: \tilde{\mathbf \Gamma}_{\text{i}}, \quad
\check{\mathbf X} : \ \stackrel{\diamondsuit} {\check{\mathbf \Gamma}}_{\text{i}} =
\hat{\mathbf X}: \stackrel{\triangle} {\hat{\mathbf \Gamma}}_{\text{i}} =
\tilde{\mathbf X} : \dot{\tilde{\mathbf \Gamma}}_{\text{i}}.
\end{equation}

\subsection{Free energy}

Suppose that the free energy is given as a sum of isotropic functions\footnote
{the first two terms of this additive split can be
motivated by the rheological model (fig. \ref{fig1}.a).}
 (cf. \cite{Lion}, \cite{Helm1})
\begin{equation}\label{freeen}
\psi=\psi(\hat{\mathbf \Gamma}_{\text{e}}, \check{\mathbf \Gamma}_{\text{ie}}, s_{\text{e}})=
\psi_{\text{el}} (\hat{\mathbf \Gamma}_{\text{e}}) + \psi_{\text{kin}} (\check{\mathbf \Gamma}_{\text{ie}}) +
\psi_{\text{iso}}(s_{\text{e}}),
\end{equation}
where the tensors $\hat{\mathbf \Gamma}_{\text{e}}$ and
$\check{\mathbf \Gamma}_{\text{ie}}$ are defined by \eqref{defgae}, $\eqref{defgaie}_2$.
Here, $\psi_{\text{el}} (\hat{\mathbf \Gamma}_{\text{e}})$ corresponds to the energy, stored due to macroscopic
elastic deformations. The "inelastic" part $\psi_{\text{kin}} (\check{\mathbf \Gamma}_{\text{ie}}) +
\psi_{\text{iso}}(s_{\text{e}})$ represents the energy, stored in the microstructure during the
viscoplastic flow due to the heterogeneity of dislocations.
The following special form of the
free energy can be used
\begin{equation}\label{spec1}
\rho_{\scriptscriptstyle \text{R}}  \psi_{\text{el}}(\hat{\mathbf \Gamma}_{\text{e}})=
\frac{k}{2}\big(\text{ln}\sqrt{\text{det} \hat{\mathbf{C}}_{\text{e}}} \big)^2+
\frac{\mu}{2} \big( \text{tr} \overline{\hat{\mathbf{C}}_{\text{e}}} - 3 \big),
\end{equation}
\begin{equation}\label{spec2}
\rho_{\scriptscriptstyle \text{R}}  \psi_{\text{kin}}(\check{\mathbf{\Gamma}}_{\text{ie}})=
\frac{c}{4}\big( \text{tr} \overline{\check{\mathbf{C}}_{\text{ie}}} - 3 \big),
\quad \rho_{\scriptscriptstyle \text{R}} \psi_{\text{iso}}(s_{\text{e}}) =\frac{\gamma}{2} \ (s_{\text{e}})^2.
\end{equation}
Here,
$k >0$, $\mu > 0$, $c \geq 0$, $\gamma \geq 0$ are material constants.
The overline $\overline{(\cdot)}$ denotes the unimodular part of a tensor
\begin{equation}\label{unim}
\overline{\mathbf{A}}:=(\det \mathbf{A})^{-1/3} \mathbf{A} \ .
\end{equation}

Denote by $\frac{\displaystyle \partial \alpha (\mathbf{A})}{\displaystyle \partial \mathbf{A}}$ the
derivative of real-valued function $\alpha$ with respect to tensor-valued argument $\mathbf{A}$ such that
\begin{equation}\label{derde}
\delta \alpha = \frac{\displaystyle \partial \alpha(\mathbf{A})}{\displaystyle \partial \mathbf{A}} : \delta \mathbf{A}.
\end{equation}

Using this notation, we introduce formally the following potential relations for stresses
$\hat{\mathbf S}$, $\check{\mathbf X}$ and for isotropic hardening $R$
\begin{equation}\label{potent}
\hat{\mathbf S}= \rho_{\scriptscriptstyle \text{R}}
\frac{\displaystyle \partial \psi_{\text{el}}(\hat{\mathbf{\Gamma}}_{\text{e}})}{\displaystyle \partial \hat{\mathbf{\Gamma}}_{\text{e}}}, \quad
\check{\mathbf X}= \rho_{\scriptscriptstyle \text{R}}
\frac{\displaystyle \partial
\psi_{\text{kin}}(\check{\mathbf{\Gamma}}_{\text{ie}})}
{\displaystyle \partial \check{\mathbf{\Gamma}}_{\text{ie}}}, \quad
R= \rho_{\scriptscriptstyle \text{R}} \frac{\displaystyle \partial \psi_{\text{iso}}(s_{\text{e}})}{\displaystyle \partial s_{\text{e}}}.
\end{equation}

\subsection{Clausius-Duhem inequality}

The Clausius-Duhem inequality imposes an additional constraint on the
material response, which states that the internal dissipation is always nonnegative.
For isothermal processes the specific internal dissipation $\delta_{\text{i}}$ takes the form (see \cite{Haupt})
\begin{equation}\label{cld}
\delta_{\text{i}} := \frac{1}{\rho_{\scriptscriptstyle \text{R}}} \tilde{\mathbf T} : \dot{\mathbf E} - \dot{\psi} \geq 0.
\end{equation}
Now, let us rewrite this expression, using relations of previous subsections.
First, we note that $\hat{\mathbf S}$ and $\check{\mathbf X}$ are isotropic functions
of $\hat{\mathbf{\Gamma}}_{\text{e}}$ and $\check{\mathbf{\Gamma}}_{\text{ie}}$, respectively.
In particular, since $\hat{\mathbf S}$ and $\hat{\mathbf{\Gamma}}_{\text{e}}$ commute, we get
\begin{equation}\label{zago}
\hat{\mathbf S} : \big( \hat{\mathbf L}_{\text{i}}^{\text{T}}
\hat{\mathbf{\Gamma}}_{\text{e}} + \hat{\mathbf{\Gamma}}_{\text{e}} \hat{\mathbf L}_{\text{i}} \big)  \stackrel{\eqref{scalprod}}{=}  2
\big(\hat{\mathbf{\Gamma}}_{\text{e}} \hat{\mathbf S} \big) : \hat{\mathbf L}_{\text{i}} \stackrel{\eqref{oldrspec}_1}{=}
2 \big(\hat{\mathbf{\Gamma}}_{\text{e}} \hat{\mathbf S} \big) : \stackrel{\triangle} {\hat{\mathbf{\Gamma}}}_{\text{i}}.
\end{equation}
Further, note that
\begin{equation}\label{zago2}
\hat{\mathbf S} : \stackrel{\triangle} {\hat{\mathbf{\Gamma}}} \ \stackrel{\eqref{addco1}}{=}
\hat{\mathbf S} : \stackrel{\triangle} {\hat{\mathbf{\Gamma}}}_{\text{e}}+
\hat{\mathbf S} : \stackrel{\triangle} {\hat{\mathbf{\Gamma}}}_{\text{i}} \stackrel{\eqref{altol}_1}{=}
\hat{\mathbf S} : \stackrel{\displaystyle \cdot} {\hat{\mathbf{\Gamma}}}_{\text{e}} + \hat{\mathbf S} :
\big( \hat{\mathbf L}_{\text{i}}^{\text{T}}
\hat{\mathbf{\Gamma}}_{\text{e}} + \hat{\mathbf{\Gamma}}_{\text{e}} \hat{\mathbf L}_{\text{i}} \big) +
\hat{\mathbf S} : \stackrel{\triangle} {\hat{\mathbf{\Gamma}}}_{\text{i}} \stackrel{\eqref{zago}}{=}
\hat{\mathbf S} : \stackrel{\displaystyle \cdot} {\hat{\mathbf{\Gamma}}}_{\text{e}} +
\big(\hat{\mathbf{C}}_{\text{e}} \hat{\mathbf S} \big) : \stackrel{\triangle} {\hat{\mathbf{\Gamma}}}_{\text{i}} .
\end{equation}
In the same way, since $\check{\mathbf X}$ and $\check{\mathbf{\Gamma}}_{\text{ie}}$ commute, we obtain
from \eqref{addmic}, $\eqref{altol}_2$, $\eqref{oldrspec}_2$
\begin{equation}\label{zago3}
\check{\mathbf X} : \stackrel{\diamondsuit} {\check{\mathbf \Gamma}}_{\text{i}} =
\check{\mathbf X} : \stackrel{\displaystyle \cdot} {\check{\mathbf{\Gamma}}}_{\text{ie}} +
\big(\check{\mathbf{C}}_{\text{ie}} \check{\mathbf X} \big) : \stackrel{\diamondsuit} {\check{\mathbf{\Gamma}}}_{\text{ii}} .
\end{equation}
Thus, we get for the stress power
\begin{multline}\label{zago4}
\hat{\mathbf S} : \stackrel{\triangle} {\hat{\mathbf{\Gamma}}} \ \stackrel{\eqref{backst2}}{=}
\hat{\mathbf S} : \stackrel{\triangle} {\hat{\mathbf{\Gamma}}} -
\hat{\mathbf X} : \stackrel{\triangle} {\hat{\mathbf{\Gamma}}}_{\text{i}} +
\check{\mathbf X} : \stackrel{\diamondsuit} {\check{\mathbf \Gamma}}_{\text{i}} \\ \stackrel{\eqref{zago2}, \eqref{zago3}}{=}
\hat{\mathbf S} : \stackrel{\displaystyle \cdot} {\hat{\mathbf{\Gamma}}}_{\text{e}} +
\big(\hat{\mathbf{C}}_{\text{e}} \hat{\mathbf S} - \hat{\mathbf X} \big) : \stackrel{\triangle} {\hat{\mathbf{\Gamma}}}_{\text{i}}
+ \check{\mathbf X} : \stackrel{\displaystyle \cdot} {\check{\mathbf{\Gamma}}}_{\text{ie}} +
\big(\check{\mathbf{C}}_{\text{ie}} \check{\mathbf X} \big) : \stackrel{\diamondsuit} {\check{\mathbf{\Gamma}}}_{\text{ii}}.
\end{multline}
Hence, the internal dissipation takes the form
\begin{multline}\label{dissi}
\delta_{\text{i}} = \frac{1}{\rho_{\scriptscriptstyle \text{R}}} \tilde{\mathbf T} : \dot{\mathbf E}  - \dot{\psi}
\stackrel{\eqref{strpow1}_2}{=}
\frac{1}{\rho_{\scriptscriptstyle \text{R}}} \hat{\mathbf S} : \stackrel{\triangle} {\hat{\mathbf{\Gamma}}} - \dot{\psi}
\stackrel{\eqref{zago4}}{=} \big( \frac{1}{\rho_{\scriptscriptstyle \text{R}}} \hat{\mathbf S} -
\frac{\displaystyle \partial \psi_{\text{el}}}{\displaystyle \partial \hat{\mathbf{\Gamma}}_{\text{e}}} \big):
\stackrel{\displaystyle \cdot} {\hat{\mathbf{\Gamma}}}_{\text{e}} +
\big(\frac{1}{\rho_{\scriptscriptstyle \text{R}}} \check{\mathbf X} - \frac{\displaystyle \partial
\psi_{\text{kin}}}{\displaystyle \partial \check{\mathbf{\Gamma}}_{\text{ie}}}\big)
: \stackrel{\displaystyle \cdot} {\check{\mathbf{\Gamma}}}_{\text{ie}} \\ +
\frac{1}{\rho_{\scriptscriptstyle \text{R}}} \big(\hat{\mathbf{C}}_{\text{e}} \hat{\mathbf S} -
\hat{\mathbf X} \big) : \stackrel{\triangle} {\hat{\mathbf{\Gamma}}}_{\text{i}} +
\frac{1}{\rho_{\scriptscriptstyle \text{R}}} \big(\check{\mathbf{C}}_{\text{ie}} \check{\mathbf X} \big)
: \stackrel{\diamondsuit} {\check{\mathbf{\Gamma}}}_{\text{ii}} - \frac{\displaystyle \partial
\psi_{\text{iso}}}{\displaystyle \partial s_{\text{e}}} \ \dot{s}_{\text{e}}.
\end{multline}
We abbreviate
\begin{equation}\label{drforc}
\hat{\mathbf \Sigma} := \hat{\mathbf{C}}_{\text{e}} \hat{\mathbf S} -
\hat{\mathbf X}, \quad
\check{\mathbf \Xi} := \check{\mathbf{C}}_{\text{ie}} \check{\mathbf X}.
\end{equation}
Finally, taking into account potential relations \eqref{potent} and definition \eqref{iso1}, we simplify
\eqref{dissi} to obtain
the Clausius-Duhem inequality in the form
\begin{equation}\label{cld2}
\rho_{\scriptscriptstyle \text{R}} \delta_{\text{i}} = \Big( \hat{\mathbf \Sigma} : \stackrel{\triangle} {\hat{\mathbf{\Gamma}}}_{\text{i}} -
R \ \dot{s} \Big) + \ \check{\mathbf \Xi} : \stackrel{\diamondsuit} {\check{\mathbf{\Gamma}}}_{\text{ii}} + \
R \ \dot{s}_{\text{d}} \geq 0.
\end{equation}

\subsection{Evolution equations}

Following the standard procedure, we formulate the evolution equations for internal variables so that
inequality \eqref{cld2} holds for arbitrary mechanical loadings (cf. \cite{Lion}, \cite{Helm1}).
\begin{equation}\label{evol}
\stackrel{\triangle}{\hat{\mathbf{\Gamma}}}_{\text{i}} := \lambda_{\text{i}}
\frac{\displaystyle \hat{\mathbf \Sigma}^{\text{D}} } {\displaystyle
\| \hat{\mathbf \Sigma}^{\text{D}} \|}, \quad
\stackrel{\diamondsuit}{\check{\mathbf{\Gamma}}}_{\text{ii}} := \lambda_{\text{i}} \ \varkappa \ \check{\mathbf \Xi}^{\text{D}},
\end{equation}
\begin{equation*}\label{evol2}
\dot{s}:= \sqrt{\frac{2}{3}} \lambda_{\text{i}}, \quad
\dot{s}_{\text{d}}:= \frac{\beta}{\gamma} \dot{s} R,
\end{equation*}
\begin{equation}\label{defi}
\| \mathbf A \| := \sqrt{ \mathbf A : \mathbf A}, \quad
\mathbf A^{\text{D}} := \mathbf A - \frac{1}{3} \text{tr}(\mathbf A) \mathbf 1,
\end{equation}
where the inelastic multiplier $\lambda_{\text{i}}$ is determined according to the Perzyna rule \cite{Perzyna2}
\begin{equation}\label{perz}
\lambda_{\text{i}}:= \frac{\displaystyle 1}{\displaystyle
\eta}\Big\langle \frac{\displaystyle 1}{\displaystyle k_0}
f \Big\rangle^{m}, \quad
f:= \| \hat{\mathbf \Sigma}^{\text{D}} \|- \sqrt{\frac{2}{3}} \big[ K + R \big], \quad
\langle x \rangle := \text{max}(x,0).
\end{equation}

Here, $\varkappa \geq 0$, $\eta \geq 0$, $m \geq 1$, $\beta \geq 0$, $K > 0$ are material parameters,
$k_0 >0$ is used to get a dimensionless term in the bracket.

Let us show that inequality \eqref{cld2} is fulfilled.
For instance, we prove that
the parenthetical term in \eqref{cld2} is nonnegative. Indeed,
\begin{equation*}\label{prove}
\Big( \hat{\mathbf \Sigma} : \stackrel{\triangle} {\hat{\mathbf{\Gamma}}}_{\text{i}} -
R \ \dot{s} \Big) = \lambda_{\text{i}} \Big( \| \hat{\mathbf \Sigma}^{\text{D}} \| - \sqrt{\frac{2}{3}} R \Big) =
\begin{cases}
    0 \quad \text{if} \ \ f < 0 \ (\lambda_{\text{i}} =0) \\
    \lambda_{\text{i}} (f + \sqrt{\frac{2}{3}} K) \quad \text{if} \ \ f > 0  \ (\lambda_{\text{i}} > 0 )
\end{cases}.
\end{equation*}
Therefore, the material model, defined in sections 2.3 and 2.5, is \emph{thermodynamically consistent}.

According to the evolution equations \eqref{evol}, the tensors
$\hat{\mathbf \Sigma}^{\text{D}}$ and $\check{\mathbf \Xi}^{\text{D}}$ are termed
\emph{the driving force for inelastic flow} and
\emph{the driving force for inelastic flow of microstructure},
respectively.
Note that both flows are incompressible. In fact,
\vspace{-3mm}
\begin{equation}\label{inco2}
 (\text{det} \mathbf F_{\text{i}})^{\displaystyle \cdot} = (\text{det} \mathbf F_{\text{i}}) \
 \text{tr} \Big(\stackrel{\triangle}{\hat{\mathbf{\Gamma}}}_{\text{i}} \Big) =0,
 \quad
  (\text{det} \mathbf F_{\text{ii}})^{\displaystyle \cdot} = (\text{det} \mathbf F_{\text{ii}}) \
 \text{tr} \Big(\stackrel{\diamondsuit}{\check{\mathbf{\Gamma}}}_{\text{ii}} \Big) =0.
\end{equation}
Under appropriate initial conditions, it follows from \eqref{inco2} that
\begin{equation}\label{inco}
\text{det} \mathbf F_{\text{i}}= \text{det} \mathbf F_{\text{ii}} =
\text{det} \mathbf C_{\text{i}}= \text{det} \mathbf C_{\text{ii}}=1.
\end{equation}
The inelastic flow takes place if the overstress $f$ is positive.
A case of rate-independent plasticity is covered by these evolution equations as
viscosity $\eta$ tends to zero.

\subsection{Transformation to the reference configuration}

A direct numerical treatment of Oldroyd derivatives in \eqref{evol}, formulated
with respect to fictitious configurations $\hat{\mathcal{K}}$ and $\check{\mathcal{K}}$,
is complicated.
In this subsection we rewrite the material model in terms of strain-like
internal variables $\mathbf C_{\text{i}}$, $\mathbf C_{\text{ii}}$, $s$, $s_{\text{d}}$ such
that the rate of $\mathbf C_{\text{i}}$, $\mathbf C_{\text{ii}}$ will be given
by the material time derivatives $\dot{\mathbf C}_{\text{i}}$, $\dot{\mathbf C}_{\text{ii}}$.
The transformation of the model includes:
\begin{itemize}
\item representation of the free energy $\psi$ through
$\mathbf C, \mathbf C_{\text{i}}, \mathbf C_{\text{ii}}, s, s_{\text{d}}$.
\item transformation of the potential relations for stresses.
\item representation of $\| \hat{\mathbf \Sigma}^{\text{D}} \|$ through $\mathbf C, \mathbf C_{\text{i}}, \mathbf C_{\text{ii}}$.
\item transformation of the evolution equations.
\end{itemize}

\subsubsection{Representation of the free energy}

Let $(J_1, J_2, J_3)$ be a full system of invariants of a second-rank tensor, defined by
\vspace{-2mm}
\begin{equation*}\label{inva}
J_1 (\mathbf A) := \text{tr} \ \mathbf A, \quad
J_2 (\mathbf A) := \frac{1}{2} \text{tr} \ \mathbf A^2, \quad
J_3 (\mathbf A) := \frac{1}{3} \text{tr} \ \mathbf A^3.
\end{equation*}
Then, using multiplicative decompositions
\eqref{mude1}, \eqref{mude2} and the property \\ $\text{tr} (\mathbf A \mathbf B) = \text{tr} (\mathbf B \mathbf A)$,
it is easily proved that
\begin{equation}\label{inva2}
J_k (\hat{\mathbf C}_{\text{e}}) = J_k (\mathbf C {\mathbf C_{\text{i}}}^{-1}), \quad
J_k (\check{\mathbf C}_{\text{ie}}) = J_k (\mathbf C_{\text{i}} {\mathbf C_{\text{ii}}}^{-1}), \quad k=1,2,3.
\end{equation}
Since $\psi_{\text{el}}$ and $\psi_{\text{kin}}$ are isotropic functions, it follows from \eqref{freeen}, \eqref{inva2} that
\begin{equation*}\label{freeen2}
\psi=\psi(\mathbf C, \mathbf C_{\text{i}}, \mathbf C_{\text{ii}}, s, s_{\text{d}})=
\psi_{\text{el}} (\mathbf C {\mathbf C_{\text{i}}}^{-1}) + \psi_{\text{kin}} (\mathbf C_{\text{i}} {\mathbf C_{\text{ii}}}^{-1}) +
\psi_{\text{iso}}(s-s_{\text{d}}).
\end{equation*}

\subsubsection{Transformation of the potential relations for stresses}

Recall that (cf. equation (9.60) in \cite{Haupt})
\begin{equation}\label{change}
\mathbf{A}^{\text{T}} \frac{\displaystyle \partial
\alpha(\mathbf{A} \mathbf{B} \mathbf{A}^{\text{T}})}{\displaystyle \partial
( \mathbf{A} \mathbf{B} \mathbf{A}^{\text{T}} )} \mathbf{A} =
\frac{\displaystyle \partial \alpha(\mathbf{A} \mathbf{B}
\mathbf{A}^{\text{T}})}{\displaystyle \partial \mathbf{B}}
\big|_{\mathbf{A} =\text{const}}.
\end{equation}
On the other hand,
\begin{equation*}\label{note1}
\tilde{\mathbf T}  \stackrel{\eqref{Kirch2}_2, \eqref{potent}_1}{=}
2 \rho_{\scriptscriptstyle \text{R}}
{\mathbf F}_{\text{i}}^{-1}
\frac{\displaystyle \partial \psi_{\text{el}}(\hat{\mathbf{C}}_{\text{e}})}
{\displaystyle \partial \hat{\mathbf{C}}_{\text{e}}} {{\mathbf F}_{\text{i}}^{-\text{T}}}, \quad
\tilde{\mathbf X}  \stackrel{\eqref{backst1}_1, \eqref{potent}_2}{=}
2 \rho_{\scriptscriptstyle \text{R}}
{\mathbf F}_{\text{ii}}^{-1}
\frac{\displaystyle \partial \psi_{\text{kin}}(\check{\mathbf{C}}_{\text{ie}})}
{\displaystyle \partial \check{\mathbf{C}}_{\text{ie}}} {{\mathbf F}_{\text{ii}}^{-\text{T}}}.
\end{equation*}
But,
\begin{equation*}\label{note2}
\hat{\mathbf{C}}_{\text{e}} = {\mathbf F}_{\text{i}}^{-\text{T}} \mathbf{C} {\mathbf F}_{\text{i}}^{-1}, \quad
\check{\mathbf{C}}_{\text{ie}} = {\mathbf F}_{\text{ii}}^{-\text{T}} \mathbf{C}_{\text{i}} {\mathbf F}_{\text{ii}}^{-1}.
\end{equation*}
Substituting ${\mathbf F}_{\text{i}}^{-\text{T}}, \mathbf{C}, \hat{\mathbf{C}}_{\text{e}}$ and
${\mathbf F}_{\text{ii}}^{-\text{T}}, \mathbf{C}_{\text{i}}, \check{\mathbf{C}}_{\text{ie}}$ for
$\mathbf{A}, \mathbf{B}, \mathbf{A} \mathbf{B} \mathbf{A}^{\text{T}}$ in \eqref{change}, we obtain
\begin{equation}\label{note3}
\tilde{\mathbf T}  =
2 \rho_{\scriptscriptstyle \text{R}}
\frac{\displaystyle \partial \psi_{\text{el}}(\mathbf C {\mathbf C_{\text{i}}}^{-1})}
{\displaystyle \partial \mathbf{C}}\big|_{\mathbf C_{\text{i}} = \text{const}}, \quad
\tilde{\mathbf X}  =
2 \rho_{\scriptscriptstyle \text{R}}
\frac{\displaystyle \partial \psi_{\text{kin}}(\mathbf C_{\text{i}} {\mathbf C_{\text{ii}}}^{-1})}
{\displaystyle \partial \mathbf C_{\text{i}}}\big|_{\mathbf C_{\text{ii}} = \text{const}}.
\end{equation}
Now let us show that $\mathbf C \tilde{\mathbf T}$ and ${\mathbf C}_{\text{i}} \tilde{\mathbf X}$
are isotropic functions of
$\mathbf C {\mathbf C_{\text{i}}}^{-1}$ and $\mathbf C_{\text{i}} {\mathbf C_{\text{ii}}}^{-1}$, respectively. Indeed,
note that
\begin{equation}\label{note5}
\frac{\displaystyle \partial \alpha(\mathbf{A} \mathbf{B}^{-1})}
{\displaystyle \partial \mathbf{A}}
\big|_{\mathbf{B} =\text{const}} = \frac{\displaystyle \partial \alpha(\mathbf{A} \mathbf{B}^{-1})}
{\displaystyle \partial (\mathbf{A} \mathbf{B}^{-1})} \mathbf{B}^{-1}, \quad
\frac{\displaystyle \partial J_k
(\mathbf A)}{\displaystyle \partial \mathbf A} = ({\mathbf A}^{\text{T}})^{k-1}, \quad k=1,2,3.
\end{equation}
Combining \eqref{note3} and \eqref{note5}, we get
\begin{equation*}\label{trans2}
\mathbf C \tilde{\mathbf T} = 2 \rho_{\scriptscriptstyle \text{R}} \sum_{k =1}^{3}
\frac{\displaystyle \partial \psi_{\text{el}}}
{\displaystyle \partial J_k(\mathbf{C} {\mathbf C_{\text{i}}}^{-1})} (\mathbf{C} {\mathbf C_{\text{i}}}^{-1})^k, \quad
\mathbf C_{\text{i}} \tilde{\mathbf X} = 2 \rho_{\scriptscriptstyle \text{R}} \sum_{k =1}^{3}
\frac{\displaystyle \partial \psi_{\text{kin}}}
{\displaystyle \partial J_k(\mathbf{C}_{\text{i}} {\mathbf C_{\text{ii}}}^{-1})} (\mathbf{C}_{\text{i}} {\mathbf C_{\text{ii}}}^{-1})^k.
\end{equation*}
If the special form \eqref{spec1}, \eqref{spec2}
for the free energy is used,
then, taking into account incompressibility relations \eqref{inco}, we obtain
\begin{equation}\label{trans4}
\tilde{\mathbf T} = k \ \text{ln}\sqrt{\text{det} (\mathbf C)} \
  \mathbf C^{-1} + \mu \ \mathbf C^{-1} (\overline{\mathbf C} \mathbf C_{\text{i}}^{-1})^{\text{D}}, \quad \quad
  \tilde{\mathbf X} = \frac{c}{2} \ \mathbf C_{\text{i}}^{-1} (\mathbf C_{\text{i}} \mathbf C_{\text{ii}}^{-1})^{\text{D}}.
\end{equation}

\subsubsection{Representation of $\| \hat{\mathbf \Sigma}^{\text{D}} \|$}
First, let us note that
\begin{equation}\label{trace}
\text{tr} \hat{\mathbf \Sigma}= \text{tr} (\hat{\mathbf C}_{\text{e}} \hat{\mathbf S} - \hat{\mathbf X}) =
\text{tr} (\mathbf C \tilde{\mathbf T} - \mathbf C_{\text{i}} \tilde{\mathbf X}).
\end{equation}
Next, we compute the inelastic contravariant and covariant pull-back of $\hat{\mathbf \Sigma}^{\text{D}}$
\begin{equation}\label{puba}
\big({{\mathbf F}_{\text{i}}^{-\text{T}}}\big)^* \hat{\mathbf \Sigma}^{\text{D}} =
{\mathbf C}_{\text{i}}^{-1} \mathbf C \tilde{\mathbf T} -
\tilde{\mathbf X} - \text{tr} \hat{\mathbf \Sigma} \ {\mathbf C}_{\text{i}}^{-1} \stackrel{\eqref{trace}}{=}
{\mathbf C}_{\text{i}}^{-1} \big( \mathbf C \tilde{\mathbf T} - \mathbf C_{\text{i}} \tilde{\mathbf X} \big)^{\text{D}},
\end{equation}
\begin{equation}\label{puba2}
\big({{\mathbf F}_{\text{i}}}\big)^* \hat{\mathbf \Sigma}^{\text{D}} =
\mathbf C \tilde{\mathbf T} {\mathbf C}_{\text{i}} -
{\mathbf C}_{\text{i}} \tilde{\mathbf X} {\mathbf C}_{\text{i}} -
\text{tr}\hat{\mathbf \Sigma} \ {\mathbf C}_{\text{i}} \stackrel{\eqref{trace}}{=}
\big( \mathbf C \tilde{\mathbf T} - \mathbf C_{\text{i}} \tilde{\mathbf X} \big)^{\text{D}} \mathbf C_{\text{i}}.
\end{equation}
Furthermore, since $\hat{\mathbf \Sigma}^{\text{D}} \in Sym$, we get
\begin{equation}\label{drif}
\| \hat{\mathbf \Sigma}^{\text{D}} \|^2=  \text{tr} (\hat{\mathbf \Sigma}^{\text{D}}  \hat{\mathbf \Sigma}^{\text{D}}) =
\text{tr} \Big\{ \Big[ \big({{\mathbf F}_{\text{i}}^{-\text{T}}}\big)^* \hat{\mathbf \Sigma}^{\text{D}} \Big]
\Big[ \big({{\mathbf F}_{\text{i}}}\big)^* \hat{\mathbf \Sigma}^{\text{D}} \Big] \Big\}.
\end{equation}
Substituting \eqref{puba} and \eqref{puba2} in  \eqref{drif} we obtain the norm of the driving force
\begin{equation}\label{nordiv}
\mathfrak{F}:=\| \hat{\mathbf \Sigma}^{\text{D}} \|=
\sqrt{\text{tr} \big[ \big( \mathbf C \tilde{\mathbf T} - \mathbf C_{\text{i}} \tilde{\mathbf X} \big)^{\text{D}} \big]^2 }.
\end{equation}

\subsubsection{Transformation of the evolution equations}

Note that
\begin{equation}\label{trace2}
\text{tr} \ \check{\mathbf \Xi}=
\text{tr} \ (\check{\mathbf{C}}_{\text{ie}} \check{\mathbf X}) =
\text{tr} (\mathbf C_{\text{i}} \tilde{\mathbf X}).
\end{equation}
Next, we compute the covariant pull-back of $\check{\mathbf \Xi}^{\text{D}}$ :
\begin{equation}\label{puba3}
\big({{\mathbf F}_{\text{ii}}}\big)^* \check{\mathbf \Xi}^{\text{D}}  \stackrel{\eqref{trace2}}{=}
\big(\mathbf C_{\text{i}} \tilde{\mathbf X} \big)^{\text{D}} \mathbf C_{\text{ii}} .
\end{equation}
Covariant Pull-back of \eqref{evol} yields
\begin{equation*}\label{puba4}
\dot{\mathbf C}_{\text{i}} = 2 \big({{\mathbf F}_{\text{i}}}\big)^*
\stackrel{\triangle}{\hat{\mathbf{\Gamma}}}_{\text{i}} \stackrel{\eqref{evol}}{=}
2 \frac{\lambda_{\text{i}}}{\mathfrak{F}} \big({{\mathbf F}_{\text{i}}}\big)^* \hat{\mathbf \Sigma}^{\text{D}}, \quad
\dot{\mathbf C}_{\text{ii}} = 2 \big({{\mathbf F}_{\text{ii}}}\big)^*
\stackrel{\diamondsuit}{\check{\mathbf{\Gamma}}}_{\text{ii}} \stackrel{\eqref{evol}}{=}
2 \lambda_{\text{i}} \varkappa \big({{\mathbf F}_{\text{ii}}}\big)^* \check{\mathbf \Xi}^{\text{D}}.
\end{equation*}
Combining this with \eqref{puba2} and \eqref{puba3}, we obtain
\begin{equation}\label{puba5}
\dot{\mathbf C}_{\text{i}} = 2 \frac{\lambda_{\text{i}}}{\mathfrak{F}}
\big( \mathbf C \tilde{\mathbf T} - \mathbf C_{\text{i}} \tilde{\mathbf X} \big)^{\text{D}} \mathbf C_{\text{i}}, \quad
\dot{\mathbf C}_{\text{ii}} =
2 \lambda_{\text{i}} \varkappa \big(\mathbf C_{\text{i}} \tilde{\mathbf X} \big)^{\text{D}} \mathbf C_{\text{ii}}.
\end{equation}

The material model is summarized in table \ref{table1}.

\begin{table}
\caption{Summary of the material model}
\begin{tabular}{|l l|}
\hline
 $\dot{\mathbf C}_{\text{i}} = 2 \frac{\displaystyle \lambda_{\text{i}}}{\displaystyle \mathfrak{F}}
 \big( \mathbf C \tilde{\mathbf T} - \mathbf C_{\text{i}} \tilde{\mathbf X} \big)^{\text{D}} \mathbf C_{\text{i}}$,
 & $\mathbf C_{\text{i}}|_{t=0} = \mathbf C_{\text{i}}^0$, $\det \mathbf C_{\text{i}}^0 =1$,  \\
 $\dot{\mathbf C}_{\text{ii}} =
2 \lambda_{\text{i}} \varkappa (\mathbf C_{\text{i}} \tilde{\mathbf X} \big)^{\text{D}} \mathbf C_{\text{ii}}$, &
$\mathbf C_{\text{ii}}|_{t=0} = \mathbf C_{\text{ii}}^0$, $\det \mathbf C_{\text{ii}}^0 =1$,  \\
$\dot{s}:= \sqrt{\frac{\displaystyle 2}{\displaystyle 3}} \lambda_{\text{i}}, \quad
\dot{s}_{\text{d}}:= \frac{\displaystyle \beta}{\displaystyle \gamma} \dot{y} R$, &
$s|_{t=0} = s^0, \ s_{\text{d}}|_{t=0} = s_{\text{d}}^0$, \\
$\tilde{\mathbf T} =
2 \rho_{\scriptscriptstyle \text{R}}
\frac{\displaystyle \partial \psi_{\text{el}}(\mathbf C {\mathbf C_{\text{i}}}^{-1})}
{\displaystyle \partial \mathbf{C}}\big|_{\mathbf C_{\text{i}} = \text{const}}$, &
$\tilde{\mathbf X} =
2 \rho_{\scriptscriptstyle \text{R}}
\frac{\displaystyle \partial \psi_{\text{kin}}(\mathbf C_{\text{i}} {\mathbf C_{\text{ii}}}^{-1})}
{\displaystyle \partial \mathbf C_{\text{i}}}\big|_{\mathbf C_{\text{ii}} = \text{const}}$, \\
$R = \gamma s_{\text{e}}$, \quad $s_{\text{e}} = s - s_{\text{d}}$, & \\
$\lambda_{\text{i}}:= \frac{\displaystyle 1}{\displaystyle
\eta}\Big\langle \frac{\displaystyle 1}{\displaystyle k_0}
f \Big\rangle^{m}, \quad
f= \mathfrak{F}- \sqrt{\frac{2}{3}} \big[ K + R \big]$, &
$\mathfrak{F}=
\sqrt{\text{tr} \big[ \big( \mathbf C \tilde{\mathbf T} - \mathbf C_{\text{i}} \tilde{\mathbf X} \big)^{\text{D}} \big]^2 }$. \\
 \hline
\end{tabular}
\label{table1}
\end{table}

\section{Integration algorithms}

The exact solution of \eqref{puba5} has under proper initial conditions
the following geometric property: $\mathbf{C}_{\text{i}}, \mathbf{C}_{\text{ii}}$ lie on the manifold
$\mathbb{M}$, defined by
\begin{equation}\label{geopro}
\mathbb{M} := \big\{ \mathbf B \in Sym: \text{det} \mathbf B =1 \big\}.
\end{equation}
Hence, system \eqref{puba5} is a \emph{system of  differential equations on the manifold} (cf. the paper \cite{Hair}).
In this section we analyse two numerical schemes, such that the numerical solution lies exactly on $\mathbb{M}$.

\subsection{Modified Euler-Backward and exponential scheme}
Consider the Cauchy problem for a system of nonlinear ordinary differential equations
\begin{equation*}\label{difur}
\dot{\mathbf A} (t) = \mathbf{f} (\mathbf A(t), t) \mathbf A (t), \quad
\mathbf A(0)=\mathbf A^0, \quad \det(\mathbf A^0)=1.
\end{equation*}
Suppose that the tensor-valued function $\mathbf{f}$ is sufficiently smooth,
$\text{tr} ( \mathbf{f} (\mathbf B, t) )=0$ and
\begin{equation}\label{fprop}
\big(\mathbf{f} (\mathbf B, t)\big)^k \mathbf B
\in Sym \quad  \forall \ \mathbf B \in Sym, \ k = 1,2,3,...
\ .
\end{equation}
Under such conditions the exact solution lies on $\mathbb{M}$.

\emph{Remark}: condition \eqref{fprop} is nontrivial, since
$\mathbf{f} (\mathbf B, t)$ is, in general, an anisotropic function of $\mathbf B$.

By ${}^n \mathbf A, {}^{n+1} \mathbf A$ denote numerical solutions respectively at $t_n$ and $t_{n+1}$, \\
$\Delta t:= t_{n+1} - t_n$.
Suppose that ${}^n \mathbf A \in \mathbb{M}$ is given.
The classical Euler-Backward method (EBM) uses the equation with respect to the unknown ${}^{n+1} \mathbf A$ :
\begin{equation}\label{Eulcl}
{}^{n+1} \mathbf A = \big[ \mathbf 1 - \Delta t \ \mathbf{f} ({}^{n+1} \mathbf A, t_{n+1}) \big]^{-1}
 \ {}^n \mathbf A.
\end{equation}
Recall that for small $\mathbf B$ \footnote{The Neumann series \eqref{Neum} converges if $\| \mathbf B \|^* < 1$.}
\begin{equation}\label{Neum}
\big[ \mathbf 1 - \mathbf B \big]^{-1} = \mathbf 1 + \mathbf B + \mathbf B^2 + \mathbf B^3 + ... \ .
\end{equation}
The exponential method (EM) is based on the equation
\begin{equation}\label{Expo}
{}^{n+1} \mathbf A = \exp\big(\displaystyle \Delta t \ \mathbf{f} ({}^{n+1} \mathbf A, t_{n+1})\big) \ {}^n \mathbf A,
\end{equation}
where the tensor exponential is given by
\begin{equation}\label{Expde}
\exp\big( \mathbf B \big) : = \mathbf 1 + \mathbf B + \frac{1}{2 !} \mathbf B^2 + \frac{1}{3 !} \mathbf B^3 + ... \ .
\end{equation}
Let us show that \emph{both methods yield a symmetric solution}.
The idea of the proof is as follows. Substituting
\eqref{Neum}
for $\big[ \mathbf 1 - \Delta t \ \mathbf{f} ({}^{n+1} \mathbf A, t_{n+1}) \big]^{-1}$ in \eqref{Eulcl}, and
\eqref{Expde} for $\exp \big(\displaystyle \Delta t \ \mathbf{f} ({}^{n+1} \mathbf A, t_{n+1})\big)$ in \eqref{Expo}, we get
for both methods
\begin{equation}\label{Sympr}
{}^{n+1} \mathbf A = \mathbf K ({}^{n+1} \mathbf A),
\end{equation}
\begin{equation}\label{Sympr0}
\mathbf K (\mathbf B) := \Big( \mathbf 1 + \Delta t \ \mathbf{f} (\mathbf B, t_{n+1})+
 \sum_{k=2}^{\infty} c_k
\big( \mathbf{f} (\mathbf B, t_{n+1}) \big)^k \Big) \ {}^n \mathbf A,
\end{equation}
with some coefficients $c_k$.
Next, let us consider an auxiliary problem
\begin{equation}\label{Sympr2}
\mathbf{A}_{\text{aux}} = \text{sym} (\mathbf K (\mathbf{A}_{\text{aux}})).
\end{equation}
Here, the symmetrization operator
$\text{sym} (\cdot) = \frac{\displaystyle 1}{\displaystyle 2} \Big( (\cdot) + (\cdot)^{\text{T}} \Big)$ is used.
Suppose $\mathbf{A}_{\text{aux}}$ is a solution of \eqref{Sympr2}.
According to properties \eqref{fprop}, since $\mathbf{A}_{\text{aux}}$
is symmetric, we obtain
\begin{equation}\label{Sympr3}
\mathbf{A}_{\text{aux}} \ {}^{n} \mathbf A^{-1} \ \mathbf{A}_{\text{aux}} \in Sym,
\quad
\mathbf K (\mathbf{A}_{\text{aux}}) \ {}^{n} \mathbf A^{-1} \ \mathbf{A}_{\text{aux}} \in Sym.
\end{equation}
Subtracting $\eqref{Sympr3}_1$ from $\eqref{Sympr3}_2$ and taking \eqref{Sympr2} into account, we get
\begin{equation}\label{Sympr4}
\text{skew} \big(\mathbf K (\mathbf{A}_{\text{aux}})\big) \ {}^{n} \mathbf A^{-1} \ \mathbf{A}_{\text{aux}} \in Sym, \quad
\text{skew} (\cdot):= (\cdot) - \text{sym}(\cdot).
\end{equation}
Here $\text{skew} (\cdot)$ stands for the skew-symmetric part of a tensor. Thus, $\eqref{Sympr4}_1$ yields
\begin{equation}\label{Sympr5}
\text{skew} \Big( \text{skew} \big(\mathbf K (\mathbf{A}_{\text{aux}})\big) \ {}^{n} \mathbf A^{-1} \
\mathbf{A}_{\text{aux}} \Big) = \mathbf 0.
\end{equation}
Since $\text{skew} \ \text{skew} (\cdot) = \text{skew} (\cdot)$, from \eqref{Sympr5} follows
\begin{equation}\label{Sympr6}
\text{skew} \big(\mathbf K (\mathbf{A}_{\text{aux}})\big) =
\text{skew} \Big( \text{skew} \big(\mathbf K (\mathbf{A}_{\text{aux}})\big) \
\big( \mathbf 1 - {}^{n} \mathbf A^{-1} \ \mathbf{A}_{\text{aux}} \big) \Big).
\end{equation}
Denote by $\| \cdot \|^{*} $ an induced norm of a tensor
\begin{equation}\label{Openo}
\| \mathbf B \|^{*} := \max_{\|\mathbf x\|_2=1} \| \mathbf B \mathbf x \|_2, \quad \|\mathbf x\|_2 := \sqrt{x^2_1+x^2_2+x^2_3}.
\end{equation}
Then,
\begin{equation}\label{Openo2}
\| \mathbf A \mathbf B \|^{*} \leq \| \mathbf A \|^{*} \| \mathbf B \|^{*}, \quad
\| \text{skew} (\mathbf A) \|^{*} \leq \| \mathbf A \|^{*}.
\end{equation}
Note also that for small $\Delta t$
\begin{equation}\label{Openo3}
\| \mathbf 1 - {}^{n} \mathbf A^{-1} \ \mathbf{A}_{\text{aux}}  \|^{*}
 \leq \frac{1}{2}.
\end{equation}
Taking the norm of both sides of \eqref{Sympr6} and using \eqref{Openo2}, \eqref{Openo3}, we get
\begin{multline*}\label{Openo4}
\| \text{skew} \big(\mathbf K (\mathbf{A}_{\text{aux}})\big) \|^{*} \stackrel{\eqref{Sympr6}, \eqref{Openo2}_2 }{\leq}
\big\|  \text{skew} \big(\mathbf K (\mathbf{A}_{\text{aux}})\big) \
\big( \mathbf 1 - {}^{n} \mathbf A^{-1} \ \mathbf{A}_{\text{aux}} \big) \big\|^{*} \\
\stackrel{\eqref{Openo2}_1}{\leq} \| \text{skew} \big(\mathbf K (\mathbf{A}_{\text{aux}})\big) \|^{*} \
\| \mathbf 1 - {}^{n} \mathbf A^{-1} \ \mathbf{A}_{\text{aux}}  \|^{*}
\stackrel{\eqref{Openo3}}{\leq}
\frac{1}{2} \| \text{skew} \big(\mathbf K (\mathbf{A}_{\text{aux}})\big) \|^{*}.
\end{multline*}
This implies that $\| \text{skew} \big(\mathbf K (\mathbf{A}_{\text{aux}})\big) \|^{*} = 0$.
Therefore, $\mathbf K (\mathbf{A}_{\text{aux}}) \in Sym$ and $\mathbf{A}_{\text{aux}}$ is a solution of
\eqref{Sympr}. In other words, equations \eqref{Sympr} and \eqref{Sympr2} are \emph{equivalent}. $\blacksquare$

This means that \emph{no modifications of
\eqref{Eulcl} and \eqref{Expo} are necessary} to ensure the symmetry of the solution ${}^{n+1} \mathbf A$.
The reader will have no difficulty in showing that \emph{the problem of symmetry does not occur} also for a system of equations
of type \eqref{Sympr}, \eqref{Sympr0}.

The following modifications of equations \eqref{Eulcl} and \eqref{Expo} leave the corresponding original solutions unchanged:
\begin{equation}\label{Eulcl2}
{}^{n+1} \mathbf A = \text{sym} \big\{ \big[ \mathbf 1 - \Delta t \ \mathbf{f} ({}^{n+1} \mathbf A, t_{n+1}) \big]^{-1}
\ {}^n \mathbf A \big\},
\end{equation}
\begin{equation}\label{Expo2}
{}^{n+1} \mathbf A = \text{sym} \big\{ \exp\big(\displaystyle \Delta t \
\mathbf{f} ({}^{n+1} \mathbf A, t_{n+1})\big) \ {}^n \mathbf A \big\}.
\end{equation}

The advantage of the exponential method based on \eqref{Expo} or \eqref{Expo2} is
that the constraint $\text{det} ({}^{n+1} \mathbf A) =1$ is exactly satisfied.
In this paper we modify the right-hand side of \eqref{Eulcl} and \eqref{Eulcl2}, using the projection
$\overline{(\cdot)} = (\det (\cdot))^{-1/3} (\cdot)$
on the group of unimodular tensors  (cf. \cite{Helm2}):
\begin{equation}\label{Eulcl3}
{}^{n+1} \mathbf A = \overline{\big[ \mathbf 1 - \Delta t \ \mathbf{f} ({}^{n+1} \mathbf A, t_{n+1}) \big]^{-1}}
 \ {}^n \mathbf A,
\end{equation}
\begin{equation}\label{Eulcl4}
{}^{n+1} \mathbf A = \overline{\text{sym} \big\{ \big[ \mathbf 1 - \Delta t \ \mathbf{f} ({}^{n+1} \mathbf A, t_{n+1}) \big]^{-1}
\ {}^n \mathbf A \big\}}.
\end{equation}
Let us remark that both \eqref{Eulcl3} and \eqref{Eulcl4} yield the same solution ${}^{n+1} \mathbf A \in \mathbb{M}$.
Thus, the modified Euler-Backward (MEBM) is formulated by \eqref{Eulcl3} or \eqref{Eulcl4}.
Further, we notice that the exponential method (EM) \eqref{Expo} is equivalent to
\begin{equation}\label{Expo4}
{}^{n+1} \mathbf A = \overline{\text{sym} \big\{ \exp\big(\displaystyle \Delta t \
\mathbf{f} ({}^{n+1} \mathbf A, t_{n+1})\big) \ {}^n \mathbf A \big\}}.
\end{equation}

\emph{Remark.}
We have a freedom in choosing between \eqref{Eulcl3} and \eqref{Eulcl4} for
MEBM. Similarly, the EM can be based either on \eqref{Expo} or \eqref{Expo4}.
In this paper we use symmetrized equations \eqref{Eulcl4}, \eqref{Expo4}.
The reason is that these two equations can be formulated with respect to six real unknowns.
At the same time equations \eqref{Eulcl3}, \eqref{Expo} are formulated with respect to nine
independent real unknowns.

\subsection{Adaptation of integration methods to the evolution equations}

Suppose that the deformation gradient ${}^{n+1} \mathbf F$ at the time $t_{n+1} = t_n + \Delta t$
is known. Further, assume that the internal variables $\mathbf C_{\text{i}}, \mathbf C_{\text{ii}}, s, s_{\text{d}}$
at the time $t_n$ are given by ${}^{n} \mathbf C_{\text{i}}, {}^{n} \mathbf C_{\text{ii}}, {}^{n} s, {}^{n} s_{\text{d}}$, respectively.
In this subsection we formulate a system of equations for finding the internal variables at the time $t_{n+1}$.

First, we adopt the modified Euler-Backward scheme \eqref{Eulcl4} and the exponential scheme \eqref{Expo4}
to the numerical integration of evolution equations \eqref{puba5}.
We stress that the right-hand sides in \eqref{puba5} satisfy requirements \eqref{fprop}.
For instance, let us analyse the evolution equation for $\mathbf C_{\text{i}}$.
Note that, since $\mathbf C \tilde{\mathbf T}$ and ${\mathbf C}_{\text{i}} \tilde{\mathbf X}$
are isotropic functions of
$\mathbf C {\mathbf C_{\text{i}}}^{-1}$ and $\mathbf C_{\text{i}} {\mathbf C_{\text{ii}}}^{-1}$,
\begin{equation*}\label{propsa1}
2 \frac{\lambda_{\text{i}}}{\mathfrak{F}} \big( \mathbf C \tilde{\mathbf T} - \mathbf C_{\text{i}} \tilde{\mathbf X} \big)^{\text{D}} =
d_1 \mathbf 1 + d_2 \mathbf C \mathbf C_{\text{i}}^{-1} + d_3 (\mathbf C \mathbf C_{\text{i}}^{-1})^2 +
d_4 \mathbf C_{\text{i}} \mathbf C_{\text{ii}}^{-1} + d_5 (\mathbf C_{\text{i}} \mathbf C_{\text{ii}}^{-1})^2,
\end{equation*}
with some suitable $d_n \in \mathbb{R}$. It remains to check that
\begin{equation*}\label{propsa2}
(\mathbf C \mathbf C_{\text{i}}^{-1})^k \mathbf C_{\text{i}} \in Sym, \
(\mathbf C_{\text{i}} \mathbf C_{\text{ii}}^{-1})^k \mathbf C_{\text{i}} \in Sym, \quad \forall \ \mathbf C_{\text{i}}, \mathbf C_{\text{ii}}
\in Sym, \ k= 1,2,3,... \ .
\end{equation*}

The evolution equations for $s, s_{\text{d}}$ are discretized by implicit Euler scheme. Further, consider
an incremental inelastic parameter
\begin{equation}\label{ineinc}
\xi := \Delta t \ {}^{n+1} \lambda_{\text{i}}.
\end{equation}
Finally, we get the following system of equations.

\begin{equation}\label{dissys1}
{}^{n+1} {\mathbf C}_{\text{i}} - \overline{\text{sym} \big( \mathbf K_{\text{i}} ({}^{n+1} {\mathbf C}, {}^{n+1} {\mathbf C}_{\text{i}},
{}^{n+1} \mathbf C_{\text{ii}},  \xi) \big)}= \mathbf 0,
\end{equation}
\begin{equation}\label{dissys2}
{}^{n+1} {\mathbf C}_{\text{ii}} - \overline{\text{sym} \big( \mathbf K_{\text{ii}} ({}^{n+1} {\mathbf C}_{\text{i}},
{}^{n+1} \mathbf C_{\text{ii}},  \xi) \big)} = \mathbf 0,
\end{equation}
\begin{equation}\label{dissys4}
\xi= \frac{\displaystyle \Delta t}{\displaystyle \eta}\Big \langle
\frac{\displaystyle {}^{n+1} f}{\displaystyle k_0} \Big\rangle^m,
\end{equation}
\begin{equation}\label{dissys3}
{}^{n+1} s = {}^n s + \sqrt{\frac{2}{3}} \xi, \quad
{}^{n+1} s_{\text{d}}  = {}^n s_{\text{d}} + \frac{\beta}{\gamma}\sqrt{\frac{2}{3}} \xi \ {}^{n+1} R,
\end{equation}
\begin{equation}\label{dissys5}
{}^{n+1} R   = \gamma ({}^{n+1} s - {}^{n+1} s_{\text{d}}), \quad
{}^{n+1} f={}^{n+1} \mathfrak{F}-\sqrt{\frac{2}{3}} ( K + {}^{n+1} R ),
\end{equation}
\begin{equation}\label{dissys6}
{}^{n+1} \mathfrak{F}=\sqrt{\text{tr}\Big[ \big({}^{n+1} \mathbf C \ {}^{n+1} \tilde{\mathbf T}
- {}^{n+1} \mathbf C_{\text{i}} \ {}^{n+1} \tilde{\mathbf X} \big)^{\text{D}}  \Big]^2},
\end{equation}
where the operators $\mathbf K_{\text{k}}, \ \text{k} \in \{\text{i},\text{ii}\}$ are defined by
\begin{equation*}\label{dissys7}
\mathbf K_{\text{k}} :=
\begin{cases}
    \big[ \mathbf 1 - \mathbf B_{\text{k}} \big]^{-1} \ {}^{n} {\mathbf C}_{k} \quad \text{if MEBM is employed } \\
    \exp \big[ \mathbf B_{\text{k}} \big] \ {}^{n} {\mathbf C}_{k} \quad \text{if EM is employed }
\end{cases},
\end{equation*}
\begin{equation}\label{dissys8}
\mathbf B_{\text{i}} ({}^{n+1} {\mathbf C}, {}^{n+1} {\mathbf C}_{\text{i}},
{}^{n+1} \mathbf C_{\text{ii}},  \xi):=
2 \frac{\xi}{{}^{n+1} \mathfrak{F}} \big({}^{n+1} \mathbf C \ {}^{n+1}\tilde{\mathbf T} -
{}^{n+1} \mathbf C_{\text{i}} \ {}^{n+1} \tilde{\mathbf X} \big)^{\text{D}},
\end{equation}
\begin{equation}\label{dissys9}
\mathbf B_{\text{ii}} ({}^{n+1} {\mathbf C}_{\text{i}},
{}^{n+1} \mathbf C_{\text{ii}},  \xi):=
2 \ \xi \ \varkappa  \big({}^{n+1} \mathbf C_{\text{i}} \ {}^{n+1} \tilde{\mathbf X} \big)^{\text{D}}.
\end{equation}
Here ${}^{n+1} \tilde{\mathbf T}, \ {}^{n+1} \tilde{\mathbf X}$ are functions of
${}^{n+1} {\mathbf C}, {}^{n+1} {\mathbf C}_{\text{i}}, {}^{n+1} {\mathbf C}_{\text{ii}}$, given by
\eqref{note3} (or by \eqref{trans4} if the special form of $\psi_{\text{el}}$, $\psi_{\text{kin}}$ is used).

\subsection{Solution strategy}
First, we exclude ${}^{n+1}s, {}^{n+1}s_{\text{d}}$ from \eqref{dissys3}, $\eqref{dissys5}_1$ (cf. \cite{Helm2}) to get
\begin{equation}\label{excl1}
{}^{n+1} R = R (\xi) :=  \frac{ {}^{\text{t}} R + \sqrt{\frac{2}{3}} \gamma \xi }{1 +  \sqrt{\frac{2}{3}} \beta \xi}, \quad
{}^{\text{t}} R: = \gamma ({}^n s - {}^n s_{\text{d}}).
\end{equation}
Next, substituting \eqref{dissys6} and \eqref{excl1} in $\eqref{dissys5}_2$, we represent ${}^{n+1} f$
as a function of ${}^{n+1} \mathbf C_{\text{i}}, {}^{n+1} \mathbf C_{\text{ii}}, \xi$.
Thus, the problem is reduced to system \eqref{dissys1}, \eqref{dissys2}, \eqref{dissys4} with respect
to ${}^{n+1} \mathbf C_{\text{i}}, {}^{n+1} \mathbf C_{\text{ii}}, \xi$.

In this paper we decompose problem \eqref{dissys1}, \eqref{dissys2}, \eqref{dissys4} as follows.
The variables ${}^{n+1} \mathbf C_{\text{i}}, {}^{n+1} \mathbf C_{\text{ii}}$ are uniquely determined by
system \eqref{dissys1}, \eqref{dissys2} with a given $\xi$. Let us denote the corresponding solution
by $\big(\mathbf C_{\text{i}} ({}^{n+1}{\mathbf C},  \xi), \mathbf C_{\text{ii}} ({}^{n+1}{\mathbf C},  \xi) \big)$.
Substituting this solution
in \eqref{dissys6}, we obtain a function $\mathfrak{F} ({}^{n+1}{\mathbf C},  \xi)$.

If $\mathfrak{F} ({}^{n+1}{\mathbf C},  0) - \sqrt{\frac{2}{3}} ( K + {}^t R ) \leq 0$, then we put $\xi :=0$,
${}^{n+1} \mathbf C_{\text{i}} :=  {}^{n} \mathbf C_{\text{i}}$,
${}^{n+1} \mathbf C_{\text{ii}} :=  {}^{n} \mathbf C_{\text{ii}}$
(no inelastic flow occurs). Otherwise,  $\xi$ is computed using
equation \eqref{dissys4}.
Substituting $\eqref{dissys5}_2$ for ${}^{n+1} f$ in \eqref{dissys4}, we
obtain two alternative forms of the incremental consistency condition:
\begin{equation}\label{consi1}
H ({}^{n+1}{\mathbf C},  \xi):= \frac{\eta \xi}{\Delta t} -
\Bigg(\frac{ \mathfrak{F} ({}^{n+1}{\mathbf C},  \xi) - \sqrt{\frac{2}{3}} ( K +  R (\xi))}{k_0} \Bigg)^m =0,
\end{equation}
\begin{equation}\label{consi2}
D ({}^{n+1}{\mathbf C},  \xi):= \Big( \frac{\eta \xi}{\Delta t} \Big)^{1/m} - \frac{
\mathfrak{F} ({}^{n+1}{\mathbf C},  \xi) - \sqrt{\frac{2}{3}} ( K +  R (\xi))}{k_0} =0.
\end{equation}

After the solution $\xi$ is found,
the values of ${}^{n+1} \mathbf C_{\text{i}}$, ${}^{n+1} \mathbf C_{\text{ii}}$ are given by
$\mathbf C_{\text{i}} ({}^{n+1}{\mathbf C},  \xi)$, $\mathbf C_{\text{ii}} ({}^{n+1}{\mathbf C},  \xi)$. Finally,
we update $s$ and $s_{\text{d}}$ using equations \eqref{dissys3}.

\emph{Remark}: Solving system \eqref{dissys4}, $\eqref{dissys5}_2$
with respect to $\xi$ with a given $\mathfrak{F}$, it is possible to
represent $\xi$ as
a function of ${}^{n+1} \mathbf C_{\text{i}}, {}^{n+1} \mathbf C_{\text{ii}}$, thus reducing the number of unknowns.
On the other hand, for small $\eta$ this approach will result in an ill-posed problem.

\subsection{Numerical implementation}

The Newton-Raphson method is used to compute
$\big(\mathbf C_{\text{i}} ({}^{n+1}{\mathbf C},  \xi), \mathbf C_{\text{ii}} ({}^{n+1}{\mathbf C},  \xi) \big)$
from \eqref{dissys1}, \eqref{dissys2}. To this end,
equations \eqref{dissys1}, \eqref{dissys2} are linearized analytically
using the coordinate-free tensor formalism proposed by Itskov (see \cite{Itskov3}, \cite{Itskov}).

\begin{figure}\centering
\psfrag{A}[m][][1][0]{$H(\xi)$}
\psfrag{B}[m][][1][0]{$\xi$}
\psfrag{C}[m][][1][0]{if $\eta=0$}
\psfrag{E}[m][][1][0]{if $\eta > 0$}
\psfrag{F}[m][][1][0]{$0$}
\psfrag{G}[m][][1][0]{slow convergence}
\psfrag{H}[m][][1][0]{not differentiable}
\psfrag{D}[m][][1][0]{$D(\xi)$}
\scalebox{0.9}{\includegraphics{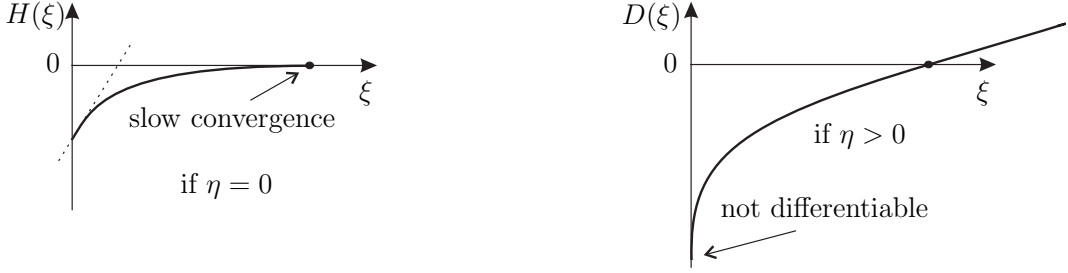}}
\caption{Finding $\xi$.\label{fig2}}
\end{figure}

Notice that the straightforward application of Newton's method to the solution of \eqref{consi1} or \eqref{consi2} is
not trivial.
Indeed, for $\eta =0, \ m > 1$, the convergence of the Newton method for \eqref{consi1} fails to be quadratic
since the first derivative is zero at the root
(see fig. \ref{fig2}). At the same time, for $\eta > 0$, the initial approximation $\xi^{(0)}=0$ can not be used to compute
the solution of \eqref{consi2}, since
the function $D(\xi)$ is not differentiable at zero (see fig. \ref{fig2}).
To overcome these difficulties, the first Newton iteration is performed using \eqref{consi1} with initial approximation $\xi^{(0)}=0$,
and the subsequent iterations are performed using \eqref{consi2}.

The derivative $\frac{\displaystyle \partial \mathfrak{F} ({}^{n+1}{\mathbf C}, \xi)}{\displaystyle \partial \xi}$,
required by the Newton method, is calculated using
the implicit differentiation of \eqref{dissys1}, \eqref{dissys2} with respect to $\xi$.
An alternative strategy is to solve \eqref{consi2} with the help of a derivative-free iteration scheme like
Pegasus method \cite{Dowell}, \cite{King}. This approach is reasonable being combined with a fixed-point iteration
for finding $\big(\mathbf C_{\text{i}} ({}^{n+1}{\mathbf C},  \xi), \mathbf C_{\text{ii}} ({}^{n+1}{\mathbf C},  \xi) \big)$,
such that no linearization of
\eqref{dissys1}, \eqref{dissys2} is required.

We implement the coordinate-free tensor formalism
to obtain an analytical expression for the consistent tangent operator
$\frac{\displaystyle \partial {}^{n+1} \tilde{\mathbf T}}{\displaystyle \partial {}^{n+1} \mathbf C}$.
Using a special product $\mathbf A \times \mathbf B$ of
two second-rank tensors and the composition $\mathbb{A} \ \mathbb{B}$ of two fourth-rank tensors
(see definitions (2.6), (2.10) in \cite{Itskov}), it follows
from \eqref{note3} that
\begin{equation*}\label{tangent}
\frac{\displaystyle \partial {}^{n+1} \tilde{\mathbf T}}{\displaystyle \partial {}^{n+1} \mathbf C} =
\frac{\displaystyle \partial \tilde{\mathbf T}
({}^{n+1} \mathbf C, {}^{n+1} \mathbf C_{\text{i}}) }{\displaystyle \partial {}^{n+1} \mathbf C} +
\frac{\displaystyle \partial \tilde{\mathbf T}
({}^{n+1} \mathbf C, {}^{n+1} \mathbf C_{\text{i}}) }{\displaystyle \partial {}^{n+1} \mathbf C_{\text{i}}}
\frac{\displaystyle \partial {}^{n+1} \mathbf C_{\text{i}} ({}^{n+1} \mathbf C) }{\displaystyle \partial {}^{n+1} \mathbf C},
\end{equation*}
\begin{equation*}\label{tangent2}
\frac{\displaystyle \partial {}^{n+1} \mathbf C_{\text{i}} ({}^{n+1} \mathbf C) }{\displaystyle \partial {}^{n+1} \mathbf C}=
\frac{\displaystyle \partial \mathbf C_{\text{i}} ({}^{n+1}{\mathbf C},  \xi) }{\displaystyle \partial {}^{n+1} \mathbf C}+
\frac{\displaystyle \partial  \mathbf C_{\text{i}} ({}^{n+1}{\mathbf C}, \xi) }{\displaystyle \partial  \xi} \times
\frac{\displaystyle \partial  \xi ({}^{n+1} \mathbf C)}{\displaystyle \partial {}^{n+1} \mathbf C},
\end{equation*}
\begin{equation*}\label{tangent3}
\frac{\displaystyle \partial  \xi ({}^{n+1} \mathbf C)}{\displaystyle \partial {}^{n+1} \mathbf C}=
- \Big(\frac{\displaystyle \partial D ({}^{n+1}{\mathbf C},  \xi)}
{\displaystyle \partial \xi}\Big)^{-1} \frac{\displaystyle \partial D
({}^{n+1} {\mathbf C},  \xi)}{\displaystyle \partial {}^{n+1} \mathbf C}.
\end{equation*}

The numerical computation of tensor exponential $\exp (\mathbf B)$
is performed
using Taylor power series expansion \eqref{Expde}.
The derivative of tensor exponential is computed
by (see \cite{ItsAks})
\begin{equation*}
\frac{\displaystyle \partial \exp (\mathbf B)}{\displaystyle \partial \mathbf B} =
\sum_{n=1}^{\infty} \frac{1}{n !} \sum_{k=0}^{n-1} \mathbf B^{n-1-k} \otimes \mathbf B^{k}.
\end{equation*}
In general, this approach fails due to
the roundoff errors, and more sophisticated techniques are required
(see, for example, \cite{Miehe}, \cite{ItsAks}, \cite{Itskov2}, \cite{Lu}).
We do not use these advanced techniques in this paper,
since in the present calculations the argument of the exponential function
is bounded. Indeed, if $\xi \leq 0.2$ then it follows from \eqref{dissys8}, \eqref{dissys9} that
\begin{equation}\label{boun}
\| \mathbf B_{\text{k}} \| \thickapprox 2 \xi \leq 0.4, \quad \text{k} \in \{\text{i},\text{ii}\}.
\end{equation}
Therefore, the roundoff errors are negligible. Moreover, under condition \eqref{boun},
the truncated power series only with few terms yield exact results up to machine precision.

\section{Numerical tests}

Now we analyse the accuracy of the integration methods presented in section 3.
Toward this end, we simulate the material behaviour under strain controlled loading.
The loading program in the time interval $t \in [0,300]$ is defined by
\begin{equation}\label{loaprog0}
\mathbf F (t) = \overline{\mathbf F^{\prime} (t)} \quad \text{or} \quad \mathbf F (t) = \mathbf F^{\prime} (t),
\end{equation}
where
\begin{equation*}\label{loaprog}
\mathbf F^{\prime} (t) :=
\begin{cases}
    (1 - t/100) \mathbf F_1  + (t/100) \mathbf F_2 \quad \text{if} \ t \in [0,100] \\
    (2 - t/100) \mathbf F_2  + (t/100-1) \mathbf F_3 \quad \text{if} \ t \in (100,200] \\
     (3 - t/100) \mathbf F_3  + (t/100-2) \mathbf F_4 \quad \text{if} \ t \in (200,300]
\end{cases},
\end{equation*}
with
\begin{equation*}\label{loaprog2}
\mathbf F_1 :=\mathbf 1, \
\mathbf F_2 := \left(
\begin{array}{ccc}
2 &  & 0 \\
0 & \displaystyle \frac{1}{\sqrt2} & 0 \\
0 & 0 & \displaystyle \frac{1}{\sqrt2}
\end{array}
\right), \
\mathbf F_3 := \left(
\begin{array}{ccc}
1 & 1  & 0 \\
0 & 1 & 0 \\
0 & 0 & 1
\end{array}
\right), \
\mathbf F_4 := \left(
\begin{array}{ccc}
\displaystyle \frac{1}{\sqrt2} &  & 0 \\
0 & 2 & 0 \\
0 & 0 & \displaystyle \frac{1}{\sqrt2}
\end{array}
\right).
\end{equation*}

\emph{Remark.} In this section we test the numerical
schemes under a variety of loading conditions, in particular,
under non-proportional loading. In this connection, the loading
programm does not have to be mechanically plausible.

The material parameters used in simulations are summarized in table \ref{table2}.

\begin{table}[h]
\caption{Material parameters}
\begin{tabular}{| l l l l |}
\hline
$k$ [MPa] & $\mu$ [MPa]  &  $c$  [MPa]  &  $\gamma$ [MPa] \\ \hline
73500 & 28200  &   3500 & 460   \\ \hline
\end{tabular} \\
\begin{tabular}{|l l l l l l|}
\hline
$K$ [MPa] & $m$ [-] & $\eta$ [$\text{s}^{-1}$] &$k_0$ [Mpa] & $\varkappa$ [$\text{MPa}^{-1}$] & $\beta$ [-]   \\ \hline
270       & 3.6     & $2 \cdot 10^6$           & 1          &  0.028                          & 5             \\ \hline
\end{tabular}
\label{table2}
\end{table}

We put the following initial conditions on the internal variables
\begin{equation}\label{Inico}
\mathbf C_{\text{i}}|_{t=0} =  \mathbf 1, \quad    \mathbf C_{\text{ii}}|_{t=0}=  \mathbf 1, \quad s|_{t=0}=0, \quad s_{\text{d}}|_{t=0}=0.
\end{equation}

Only the uniform time stepping is used in this paper.
The numerical solution obtained with extremely small time step ($\Delta t = 0.01 \text{s}$) will be named the \emph{exact solution}.

\begin{figure}\centering
\psfrag{A}[m][][1][0]{$T_{11}$}
\psfrag{B}[m][][1][0]{$T_{12}$}
\psfrag{C}[m][][1][0]{$t$}
\psfrag{D}[m][][1][0]{$\Delta t = 10 \text{s}$}
\psfrag{E}[m][][1][0]{$\Delta t = 5 \text{s}$}
\psfrag{F}[m][][1][0]{$\Delta t = 2.5 \text{s}$}
\psfrag{G}[m][][1][0]{MEBM}
\psfrag{H}[m][][1][0]{EM}
\psfrag{K}[m][][1][0]{exact solution}
\scalebox{0.9}{\includegraphics{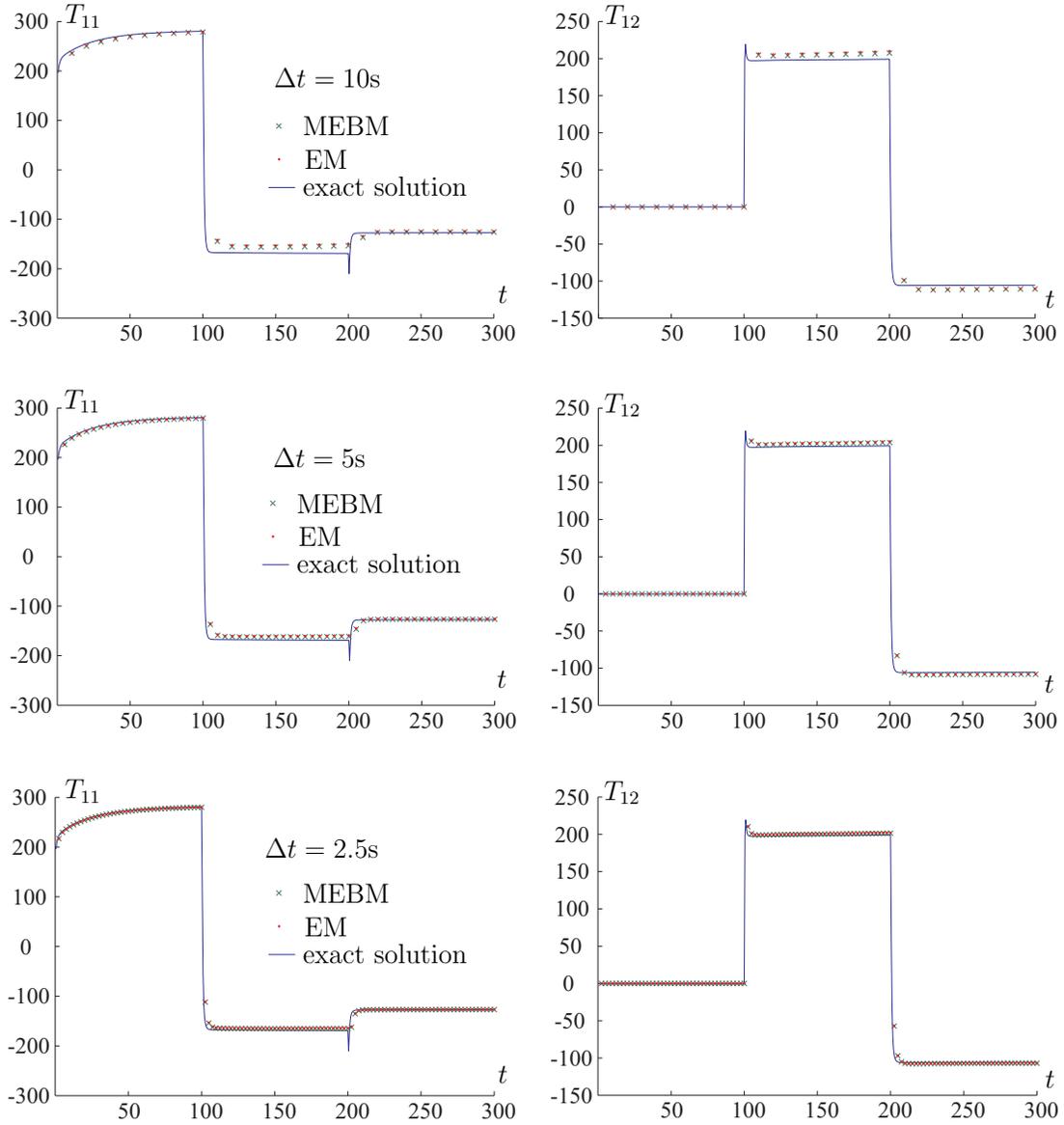}}
\caption{Accuracy test with small elastic strains (use $\eqref{loaprog0}_1$). \label{fig3}}
\end{figure}

\begin{figure}\centering
\psfrag{A}[m][][1][0]{$T_{11}$}
\psfrag{B}[m][][1][0]{$T_{12}$}
\psfrag{C}[m][][1][0]{$t$}
\psfrag{D}[m][][1][0]{$\Delta t = 10 \text{s}$}
\psfrag{E}[m][][1][0]{$\Delta t = 5 \text{s}$}
\psfrag{F}[m][][1][0]{$\Delta t = 2.5 \text{s}$}
\psfrag{G}[m][][1][0]{MEBM}
\psfrag{H}[m][][1][0]{EM}
\psfrag{K}[m][][1][0]{exact}
\scalebox{0.9}{\includegraphics{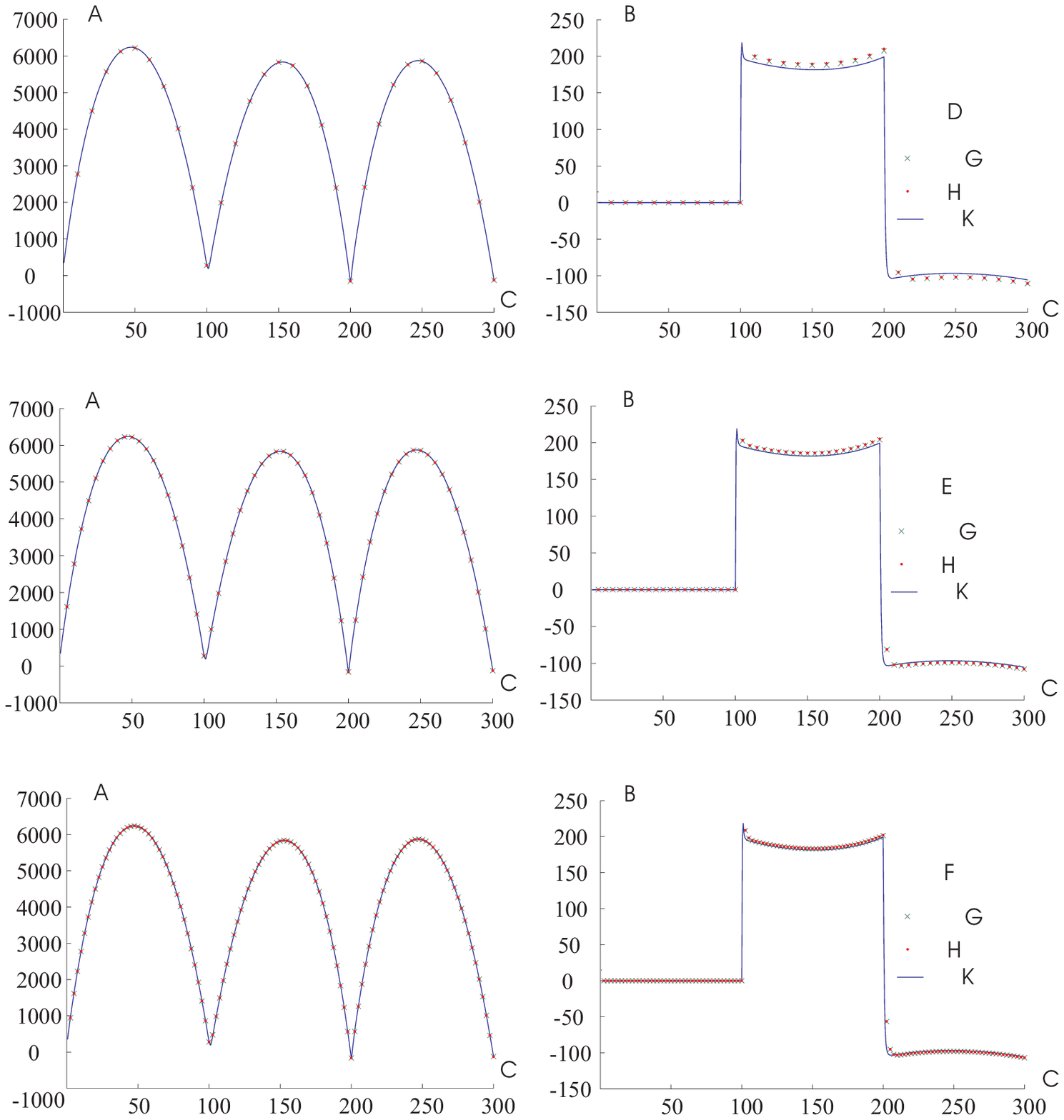}}
\caption{Accuracy test with finite elastic strains (use $\eqref{loaprog0}_2$). \label{fig4}}
\end{figure}

The coordinates of Cauchy stress tensor $\mathbf T$ for loadings $\eqref{loaprog0}_1$ and $\eqref{loaprog0}_2$ are
plotted respectively in figures \ref{fig3} and \ref{fig4}.
Note that $\det (\mathbf F) \equiv 1 $ if $\eqref{loaprog0}_1$ is used, and no hydrostatic stress occurs. On the other hand,
relation $\eqref{loaprog0}_2$ results in a large hydrostatic stress and a \emph{finite elastic bulk strain}.

The numerical simulation shows that both MEBM and EM have a similar error.
Both methods produce slightly different results for $\Delta t = 10 \ \text{s}$ when the inelastic increment $\xi$
ranges up to about 17\%.

\section{Characterization of the material model}

We investigate qualitatively the material response,
predicted by the material model.
The numerical computations
simulate basic material testing experiments.
Material parameters
from table \ref{table2} and initial conditions \eqref{Inico} are used in this section.

\subsection{Uniaxial testing}
\begin{figure}\centering
\psfrag{A}[m][][1][0]{\footnotesize $\dot{\varepsilon}=10^{-6}$/s}
\psfrag{B}[m][][1][0]{\footnotesize $\dot{\varepsilon}=0.01$/s}
\psfrag{C}[m][][1][0]{\footnotesize $\dot{\varepsilon}=0.1$/s}
\psfrag{D}[m][][1][0]{\footnotesize $\dot{\varepsilon}=1$/s}
\psfrag{E}[m][][1][0]{\footnotesize $\dot{\varepsilon}=10$/s}
\psfrag{F}[m][][1][0]{\footnotesize $\dot{\varepsilon}=10^{-6}$/s}
\psfrag{X}[m][][1][0]{\footnotesize with relaxation}
\psfrag{Z}[m][][1][0]{\footnotesize $\dot{s}=27$ MPa/s}
\psfrag{G}[m][][1][0]{\footnotesize with creep}
\psfrag{H}[m][][1][0]{\footnotesize $\dot{\varepsilon}=0.1$/s}
\psfrag{K}[m][][1][0]{a}
\psfrag{L}[m][][1][0]{b}
\psfrag{M}[m][][1][0]{c}
\psfrag{N}[m][][1][0]{$\varepsilon$}
\psfrag{p}[m][][1][0]{$\sigma$}
\scalebox{0.95}{\includegraphics{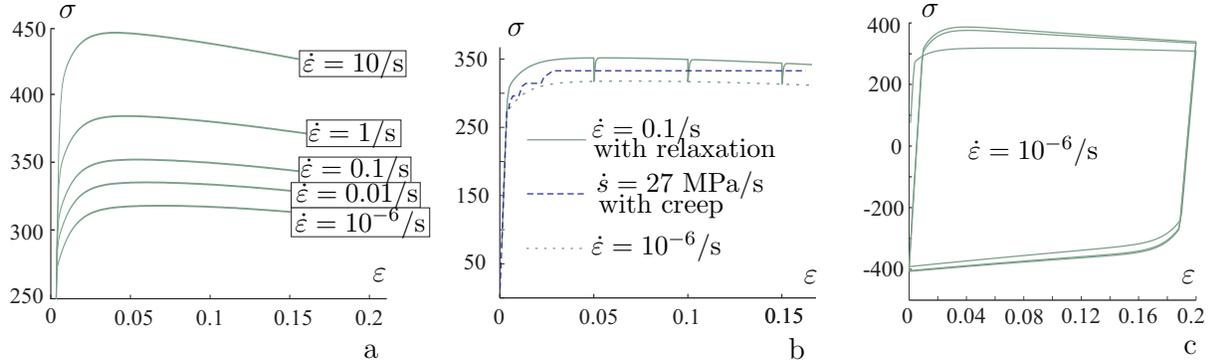}}
\caption{Uniaxial testing: monotonic loading (a), relaxation and creep (b), and cyclic
loading (c). \label{fig5}}
\end{figure}

For uniaxial test we put
\begin{equation}\label{uniax}
\mathbf F = \left(
\begin{array}{ccc}
1 + \varepsilon & 0  & 0 \\
0 & \alpha & 0 \\
0 & 0 & \alpha
\end{array}
\right), \quad T_{22}= T_{33}=0.
\end{equation}
The unknown $\alpha$ is determined using   $\eqref{uniax}_2$.
The technical stress
$\sigma:= \frac{\displaystyle A}{\displaystyle A_0} T_{11}$
is plotted in figure \ref{fig5}.a for various strain rates
$\dot{\varepsilon}$. Here $A$, $A_0$ denote the
current and initial cross sections, respectively.
Although the material response is stable, the stress reduction
is observed after the peak load in uniaxial monotonic test.
The reason is
the reduction of the cross-section.
The equilibrium curve can be reached both by
relaxation and creep (figure \ref{fig5}.b).
In the simulation presented in figure \ref{fig5}.b each
relaxation period lasts for 10 seconds. The creep time is 20 seconds.
Therefore, the numerical experiment shows that it takes
longer to reach equilibrium curve
in the creep process than in the relaxation
process. Finally, as indicated by
the strain-controlled cyclic test (figure \ref{fig5}.c), the
saturation is achieved after the isotropic hardening is accomplished.

\subsection{Torsion testing}

The cyclic torsion testing has much potential for providing
information about the nonlinear hardening phenomena (see \cite{Goerke}).
To simulate the torsion of constrained thin-walled tube we put
\begin{equation}\label{tors}
\mathbf F = \left(
\begin{array}{ccc}
1 & \phi  & 0 \\
0 & 1 & 0 \\
0 & 0 & \alpha
\end{array}
\right), \quad T_{33}=0.
\end{equation}
We consider a strain controlled torsion test with a given $\phi$, $| \dot{\phi} | = 0.01/s$.
The unknown $\alpha$ is determined using  $\eqref{tors}_2$. Denote by
$\sigma := \frac{\displaystyle A}{\displaystyle A_0} T_{22}$ and
$\tau := \frac{\displaystyle A}{\displaystyle A_0} T_{12}$ the
axial and shear stresses, respectively.

\begin{figure}\centering
\psfrag{A}[m][][1][0]{$\tau$}
\psfrag{B}[m][][1][0]{$\sigma$}
\psfrag{C}[m][][1][0]{$\phi$}
\psfrag{D}[m][][1][0]{monotonic loading}
\psfrag{E}[m][][1][0]{cyclic loading}
\psfrag{F}[m][][1][0]{a}
\psfrag{G}[m][][1][0]{b}
\psfrag{H}[m][][1][0]{c}
\psfrag{K}[m][][1][0]{d}
\scalebox{0.9}{\includegraphics{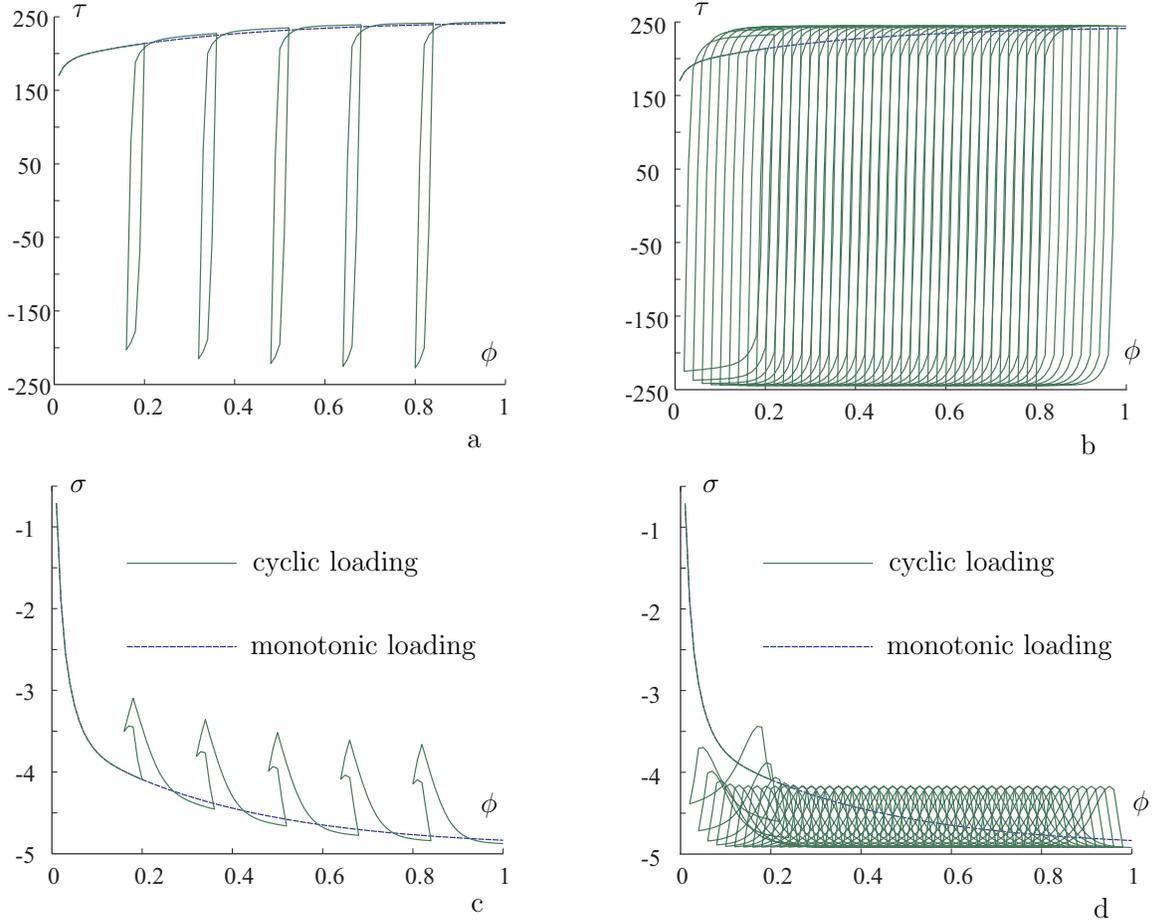}}
\caption{Torsion testing (material parameters from table \ref{table2}): shear stress (a), (b), and axial stress (c), (d).
\label{fig6}}
\end{figure}

The axial stress $\sigma$ is exactly zero in geometric linear theory. But,
in the case of finite strains, so-called second order effects can appear,
leading to nonzero axial stress.
For instance, the Poynting effect (see \cite{Andron}) is observed during the torsion
of cylindrical samples made of aluminium alloy. This effect consists in axial compression of
constrained samples or axial elongation of unconstrained samples (see \cite{Poynt}).
The axial and shear stresses are plotted in figure \ref{fig6} for different forming increments.
As may be seen from the figure, the Poynting effect is predicted by the material model.

\begin{figure}\centering
\psfrag{A}[m][][1][0]{$\tau$}
\psfrag{B}[m][][1][0]{$\sigma$}
\psfrag{C}[m][][1][0]{$\phi$}
\psfrag{F}[m][][1][0]{a}
\psfrag{G}[m][][1][0]{b}
\psfrag{H}[m][][1][0]{c}
\psfrag{K}[m][][1][0]{d}
\scalebox{0.9}{\includegraphics{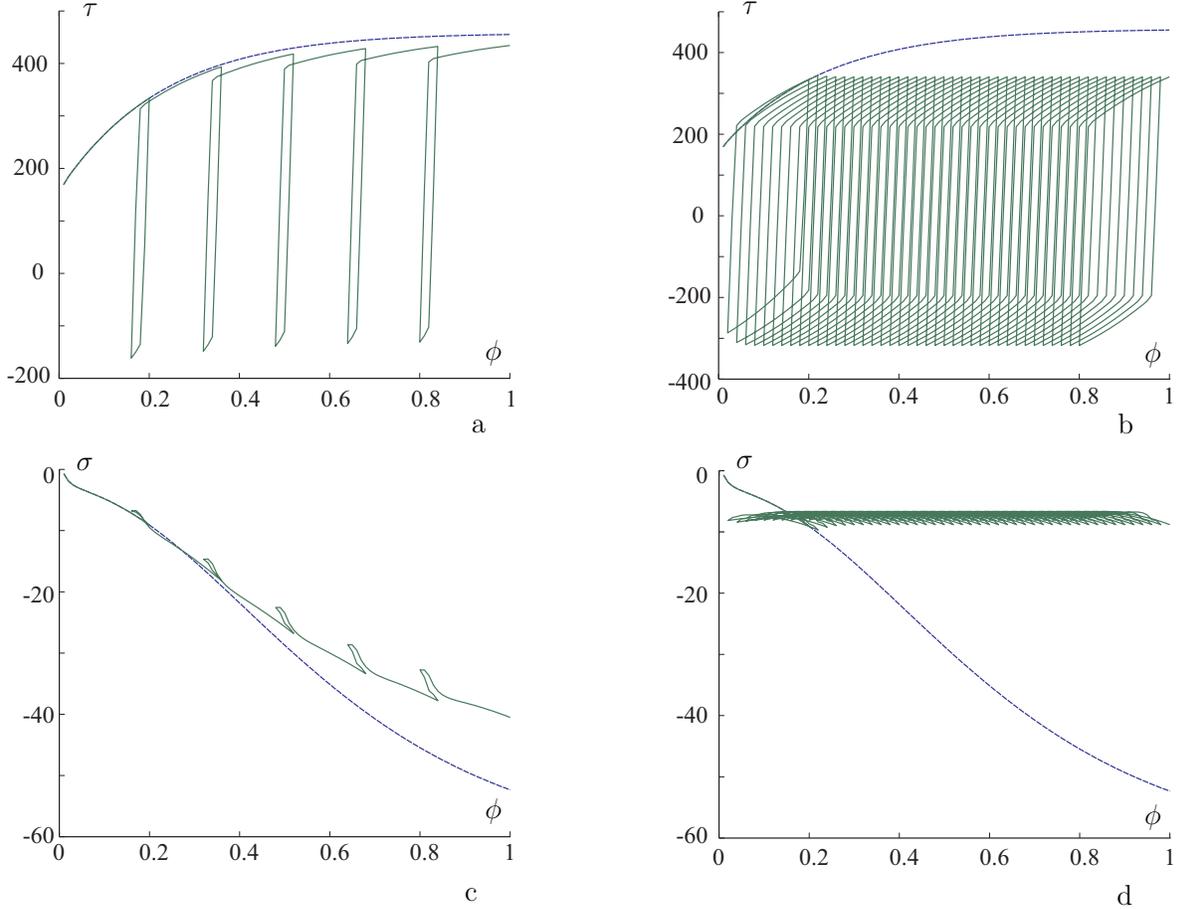}}
\caption{Torsion testing (modified material parameters): shear stress (a), (b), and axial stress (c), (d).
\label{fig7}}
\end{figure}

Next, as indicated by figure \ref{fig6}, the maximal stresses are influenced by the forming increment. For
smaller forming increment the isotropic hardening is
accomplished on the early stage of the forming process, thus leading to higher maximal stresses.
On the other hand, if the kinematic hardening
is not accomplished within one forming increment, a somewhat different material response is possible.
For instance, the maximal stresses under the cyclic loading can be smaller
than the stresses under the monotonic loading.
In order to demonstrate this effect, we perform the numerical simulation (see figure \ref{fig7})
 with modified hardening parameters:
$\varkappa := 0.0035 \ \text{MPa}^{-1}$, $c := 1500 \ \text{MPa}$, $\beta := 10$, $\gamma := 1800 \ \text{MPa}$.
A similar effect was reported in \cite{Meyer} for 20MoCr24 steel alloy.

\section{Discussion}

The classical material model of viscoplasticity is modified in a thermodynamically
consistent manner
to incorporate finite elastic and inelastic strains. The model takes
rate-dependence (relaxation, creep) and
hysteresis effects (nonlinear kinematic and isotropic hardening) into account.

Although the material response is anisotropic, the symmetry of ${}^{n+1} \mathbf C_{\text{i}}$ and
${}^{n+1} \mathbf C_{\text{ii}}$ is a priori preserved
by EBM, MEBM and EM. It is shown that \emph{no
symmetrization procedure is necessary}. Moreover,
any symmetrization should leave the corresponding solutions unchanged.
Both MEBM and EM have the advantage that the
inelastic incompressibility constraint is exactly satisfied.

Under special assumptions on the potential
functions $\psi_{\text{el}}$ and $\psi_{\text{kin}}$ it may be beneficial
to optimize the solution procedure of system
\eqref{dissys1} --- \eqref{dissys9}, thus reducing the computational effort.
On the other hand,
the most important properties of any stress algorithm are
stability, accuracy, robustness and universality.
The computational effort, required for the evaluation of
stresses and tangential operator, is negligible
in comparison with the costs of solving the global linearized
system of equations within the Newton-Raphson iterative
procedure.

The material model reproduces qualitatively the experimental results \cite{Ha1}, \cite{Ha2} for
aluminium alloy processed by ECA-pressing.
For more detailed modeling of kinematic hardening it is possible to
introduce several Armstrong-Frederick terms,
using series of multiplicative decompositions of type \eqref{mude2}.
To complete the phenomenological description of the material,
a proper parameter identification is required.

\section*{Acknowledgements}

This research was supported by German National Science Foundation (DFG) within the
collaborative research center SFB 692 "High-strength aluminium based light weight materials for
reliable components". The authors are grateful to Dr. D. Helm and Dr. P. Neff for fruitful
discussions.

\end{document}